\documentclass[11pt]{amsart} 
\usepackage[lmargin=1in,rmargin=1in,tmargin=1in,bmargin=1in]{geometry}
\usepackage[ps,all,arc,rotate]{xy}
\usepackage{graphicx, float, epstopdf}
\usepackage{hyperref, color}
\usepackage{centernot}
\usepackage{fancyhdr}
\usepackage[utf8]{inputenc}
\usepackage{amsfonts,amssymb,amsmath,amsthm}
\usepackage[usenames,dvipsnames]{xcolor}
\usepackage{datetime}
\numberwithin{equation}{section}
\numberwithin{figure}{section}
\allowdisplaybreaks[4]                        
 
\newtheorem{lemma}{Lemma}[section]
\newtheorem{theorem}{Theorem}[section]

\theoremstyle{definition}

\newtheorem{remark}{Remark}[section]
\newcommand{\R}{\mathbb{R}}

\newcommand{\N}{\mathbb{N}}

\newcommand{\e}{\operatorname{e}}
\newcommand{\real}{\operatorname{Re}}

\newcommand{\sumtwo}{\operatorname*{\sum\sum}}

\newcommand{\modu}{\operatorname{mod}}

\begin{document}
\title[Perturbed moments and a longer mollifier for critical zeros of $\zeta$]{Perturbed moments \\ and a longer mollifier for critical zeros of $\zeta$}
\author{Kyle Pratt}
\address{Department of Mathematics, University of Illinois, 1409 West Green Street, Urbana, IL 61801, United States}
\email{kpratt4@illinois.edu} 
\author{Nicolas Robles}
\address{Department of Mathematics, University of Illinois, 1409 West Green Street, Urbana, IL 61801, United States}
\email{nirobles@illinois.edu}
\address{Wolfram Research Inc, 100 Trade Center Dr, Champaign, IL 61820, USA}
\email{nicolasr@wolfram.com}
\maketitle
\begin{abstract}
Let $A(s)$ be a general Dirichlet polynomial and $\Phi$ be a smooth function supported in $[1,2]$ with mild bounds on its derivatives. New main terms for the integral $I(\alpha,\beta)=\int_{\mathbb{R}} \zeta(\frac{1}{2}+\alpha+it)\zeta(\frac{1}{2}+\beta+it)|A(\frac{1}{2}+it)|^2 \Phi(\frac{t}{T})dt$ are given. For the error term, we show that the length of the Feng mollifier can be increased from $\theta < \frac{17}{33}$ to $\theta < \frac{6}{11}$ by decomposing the error into Type I and Type II sums and then studying the resulting sums of Kloosterman sums. As an application, we slightly increase the proportion of zeros of $\zeta(s)$ on the critical line.
\end{abstract}
\section{Introduction}
\subsection{Background and motivation}
Let $A(s)$ be the Dirichlet polynomial
\begin{align} \label{defAmisc}
A(s) := \sum\limits_{n \le N} \frac{a_n}{n^s}, \quad \textnormal{with} \quad a_n \ll_\varepsilon n^\varepsilon, \quad N := T^\theta, \quad  \textnormal{and} \quad \theta<1. 
\end{align}
Research on the twisted second moment
\begin{align} \label{Iintegral}
I := \int_{-\infty}^\infty |\zeta(\tfrac{1}{2}+it)|^2 |A(\tfrac{1}{2}+it)|^2 \Phi \bigg(\frac{t}{T}\bigg)dt,
\end{align}
with $\Phi(x)$ a smooth function supported in $[1,2]$ and with derivatives satisfying $\Phi^{(j)}(x) \ll_j \log^j T$, has been well studied in the literature of the Riemann zeta-function, see e.g. \cite{bchb, bcr, conrey83a, conrey89, levinson}.  The applications of $I$ are very deep, as one may use asymptotic estimates for $I$ to make sense of the distribution of values of $L$-functions, the location of their critical zeros, as well as upper and lower bounds for the size of $L$-functions (see, among many examples, \cite{conreygg86, cis1, cis2, feng, krz01, makslimitations, rrz01}).

As often stressed, one key aspect to obtaining good results is to make sure that $\theta$ be as large as possible. One notorious example of such a benefit is that the larger $\theta$ is, the larger the proportion of zeros of $\zeta(s)$ on the critical line becomes, up to certain limitations. For example, it is known that if one could take $\theta = 1 - \varepsilon$, then the Lindel\"of hypothesis follows (see e.g. \cite{bcr}). Moreover, as shown in \cite{bg}, if one could take $\theta = \infty$ in the Conrey-Levinson mollifier (see below), then the Riemann hypothesis would follow.

For values of $\theta$ such that $\theta < \frac{1}{2}$, the literature goes back, at least, to Levinson \cite{levinson}. Indeed, taking $\theta<\frac{1}{2}$ is not at all taxing, and it is powerful enough to show that at least a third of non-trivial zeros of $\zeta(s)$ are on the critical line. Refinements on the value of $\theta$ due to Conrey \cite{conrey89} have increased that percentage to $40.88\%$. Adding, or refining, the structure of the coefficients $a_n$ of the Dirichlet polynomial in \eqref{defAmisc} also leads to improved values of the above mentioned proportion, \cite{bcy, conrey83a, feng, krz01, krz02, rrz01}.

One of the first systematic results on $I$ was produced by Balasubramanian, Conrey and Heath-Brown in \cite{bchb}. For $\theta<\frac{1}{2}$, they showed that
\begin{align} \label{resultbchb}
I = T \sum_{d, e \le N} \frac{a_d \overline{a}_e}{[d,e]} \bigg(\log \frac{T(d,e)^2}{2 \pi de} + 2 C_0 + \log 4 -1 \bigg) + o(T),
\end{align}
where $(d,e)$ and $[d,e]$ are the gcd and the lcm of $d$ and $e$, respectively. Here $C_0$ is Euler's constant. 

When $A(s)$ is a mollifier, a loose term to indicate that $A(s)$ approximately replicates the behavior of $\zeta(s)^{-1}$, they showed that one could increase $\theta$ from $\theta<\frac{1}{2}$ to $\theta<\frac{9}{17}=\frac{1}{2}+\frac{1}{34}$. This improvement allowed them to show that at least $38\%$ of the zeros are on the critical line. 

An important and subtle change of behavior takes place when $\theta>\frac{1}{2}$. This is because when $\theta<\frac{1}{2}$ only the ``diagonal'' terms contribute to the main term, and the rest is absorbed in the error; for $\theta>\frac{1}{2}$ there is a nontrivial contribution from the ``off-diagonal'' terms.

Bettin, Chandee, and Radziwi\l\l \ \cite{bcr} succeeded in breaking the $\frac{1}{2}$ barrier for an arbitrary Dirichlet polynomial. They showed that if $\theta<\frac{1}{2}+\delta$ with $\delta=\frac{1}{66}$, then
\begin{align} \label{resultbcr}
I = \sum_{d, e \le N} \frac{a_d \overline{a}_e}{[d,e]} \int_{-\infty}^\infty \bigg(\log \frac{t(d,e)^2}{2 \pi de} + 2 C_0 \bigg) \Phi \bigg(\frac{t}{T}\bigg)dt + O(T^{\frac{3}{20}} N^{\frac{33}{20}}+T^{\frac{1}{3}+\varepsilon}).
\end{align}
Observe that $\frac{1}{2} + \delta = \frac{17}{33}$.

The key to the error term in \eqref{resultbchb} was a result of Bettin and Chandee \cite{bc} on bounds of generic trilinear Kloosterman sums, see $\mathsection$2.\\ 

Let us now move on to the details of this paper. Suppose that $\mu(n)$ and $\Lambda(n)$ denote the usual M\"{o}bius and von Mangoldt functions, respectively. Using the insights and methods of \cite{bcr}, we study the more general twisted second moment 
\begin{align} \label{generalI}
I(\alpha ,\beta ) := \int_{ - \infty }^\infty  \zeta (\tfrac{1}{2} + \alpha  + it)\zeta (\tfrac{1}{2} + \beta  - it)|A(\tfrac{1}{2} + it)|^2\Phi \bigg( \frac{t}{T} \bigg)dt,
\end{align}
where $\alpha,\beta \ll (\log T)^{-1}$. We consider three cases for the coefficients $a_n$, namely
\begin{equation} \label{3cases}
a_n = \begin{cases}
O(n^\varepsilon), & \mbox{ for } \theta<\frac{17}{33}, \\
\mu^2(n) (\mu * \Lambda^{*k})(n)f(n), & \mbox{ with $f(n) \in \mathcal{F}$ and for } \theta<\frac{6}{11}, \\
\mu(n)f(n), & \mbox{ with $f(n) \in \mathcal{F}$ and for } \theta<\frac{4}{7}. 
\end{cases}
\end{equation}
Here $\mathcal{F}$ denotes the class of smooth functions given by
\[
f(n) = P\bigg( \frac{\log (N/n)}{\log N} \bigg),
\]
where $P(x)$ is a polynomial.

The third case $\mu(n)f(n)$ of \eqref{3cases} was studied by Conrey in 1989, and we call this third case the ``Levinson-Conrey,'' or simply ``Conrey,'' mollifier. The main innovations of this paper lie in studying the second case of \eqref{3cases}, and extending the range of $\theta$ for which one may prove an asymptotic formula. The second case  is colloquially known as the ``Feng'' mollifier, since it was first exploited in \cite{feng, krz01}.

The reason behind the choice of \eqref{generalI} needs to be explained. As it will be elaborated in $\mathsection$5, the presence of the terms $\alpha$ and $\beta$ is due to the fact that in order to compute the percentage of zeros on the critical line, one first computes $I(\alpha,\beta)$ and then sets $\alpha=\beta=-R/\log T$, where $R$ is a bounded positive real number of our choice (to be optimized). Therefore, the integral $I$ in \eqref{Iintegral} needs to be generalized to accommodate the variables $\alpha$ and $\beta$.

The strength of the result appearing in \cite{bcr} is the generality of $a_n$. However, often times applications of $I$ or $I(\alpha,\beta)$ allow for specialization of the shape of $a_n$. In fact, the mollifier encountered earlier to produce percentages of zeros requires that $a_n$ should be close to the M\"{o}bius function $\mu(n)$. The precise bonus coming from the shape $\mu^2(n) (\mu * \Lambda^{*k})(n)f(n)$ will be explained in $\mathsection$5. 

\subsection{Main result}
Set $L=\log T$. We are now in a position to state the results of the paper.
\begin{theorem} \label{theorem1}
Let $\alpha, \beta \ll L^{-1}$. Then one has
\begin{align}
  I(\alpha ,\beta ) = \sum\limits_{d,e \le N}  \frac{a_d \overline{a}_e}{[d,e]}\frac{(d,e)^{\alpha  + \beta }}{d^\alpha e^\beta }\int_{ - \infty }^\infty   \bigg( \zeta (1 + \alpha  + \beta ) + \zeta (1 - \alpha  - \beta )\bigg( \frac{2\pi de}{t(d,e)^2} \bigg)^{\alpha  + \beta } \bigg)\Phi\bigg(\frac{t}{T}\bigg)dt  + O(\mathcal{E}), \nonumber
\end{align}
with $\mathcal{E}$ given by the following cases
\begin{equation}
\mathcal{E} = \begin{cases}
T^{\frac{3}{20}} N^{\frac{33}{20}}+N^{1/2}T^{\frac{1}{2}+\varepsilon}, & \mbox{ if } \quad a_n \ll n^\varepsilon, \\
T^{\varepsilon}(N^{\frac{11}{6}}+N^{\frac{11}{12}}T^{\frac{1}{2}}), & \mbox{ if } \quad a_n =\mu^2(n) (\mu * \Lambda^{*k})(n)f(n), \\
T^{\varepsilon}(N^{\frac{7}{4}}+N^{\frac{7}{8}}T^{\frac{1}{2}}), & \mbox{ if } \quad a_n =\mu(n)f(n), \nonumber
\end{cases}
\end{equation}
with $f \in \mathcal{F}$. 
\end{theorem}
\begin{remark}
Let $I_M(\alpha,\beta)$ denote the main term of Theorem \ref{theorem1}. By the use of the Laurent series of $\zeta(1+s)$ around $s=0$ we have
\[
\zeta(1+s) = \frac{1}{s} + C_0 + O(s).
\]
This implies that
\[
\zeta(1+s)+\zeta(1-s)y^s = 2C_0 - \log y + O(s).
\]
Therefore when $\alpha, \beta \to 0$ we obtain
\[
\mathop {\lim }\limits_{\substack{\alpha  \to 0 \\ \beta \to 0}} I_M(\alpha ,\beta ) =  \sum\limits_{d,e \le N}  \frac{a_d \overline{a}_e}{[d,e]}\int_{ - \infty }^\infty   \bigg( \log \frac{t(d,e)^2}{2\pi d e} + 2C_0  \bigg)\Phi \bigg( \frac{t}{T} \bigg)dt =: I_M,
\]
which is the main term from \eqref{resultbcr}.
\end{remark}
\begin{remark}
Colloquially, this means that the length of the Feng mollifier can be ``pushed'' from $\theta<\frac{17}{33}$ to $\theta<\frac{6}{11}$. We succeed by exploiting more of the structure of the mollifier, rather than relying on estimates for generic forms in Kloosterman sums.
\end{remark}
\subsection{Additional results}
Sometimes it is useful to consider moment integrals where there is a cross product of two different Dirichlet polynomials. More specifically, suppose we have two Dirichlet polynomials of the form
\begin{align} \label{2dirichletseries}
A(s) := \sum_{n \le N} \frac{a_n}{n^s}, \quad B(s) := \sum_{k \le K} \frac{b_k}{k^s}, \quad \textnormal{where} \quad a_n \ll n^\varepsilon, \quad \textnormal{and} \quad b_k \ll k^\varepsilon.
\end{align}
In this case $N := T^{\theta_1}$ and $K := T^{\theta_2}$. Now we focus on integrals of the form
\[
\Upsilon(\alpha,\beta) := \int_{-\infty}^\infty \zeta(\tfrac{1}{2}+\alpha+it)\zeta(\tfrac{1}{2}+\beta-it) A \overline {B}(\tfrac{1}{2}+it) \Phi \bigg(\frac{t}{T}\bigg)dt.
\]
In this scenario, we have the following result.
\begin{theorem} \label{theorem2}
Let $\alpha, \beta \ll L^{-1}$. Then one has
\begin{align}
  \Upsilon(\alpha ,\beta ) = \sum\limits_{\substack{n \leq N \\ k \leq K}}  \frac{a_n \overline{b}_k}{[n,k]}\frac{(n,k)^{\alpha  + \beta }}{n^\alpha k^\beta }\int_{ - \infty }^\infty   \bigg( \zeta (1 + \alpha  + \beta ) + \zeta (1 - \alpha  - \beta )\bigg( \frac{2\pi nk}{t(n,k)^2} \bigg)^{\alpha  + \beta } \bigg)\Phi\bigg(\frac{t}{T}\bigg)dt  + O(\mathcal{E}), \nonumber
\end{align}
with $\mathcal{E}$ given by
\begin{equation}
\mathcal{E} = \begin{cases}
T^{\frac{1}{2}+\varepsilon} (NK)^{\frac{1}{4}} + T^{\frac{3}{20}+\varepsilon} (NK)^{\frac{7}{10}}(N+K)^{\frac{1}{4}} \\ \ \ \ \ +T^\varepsilon (NK)^{\frac{7}{8}}(N+K)^{\frac{1}{8}},  &\text{ if } \ \ \ a_n, b_n \ll n^\varepsilon, \\
T^{\varepsilon}(NK)^{\frac{3}{4}}(N+K)^{\frac{1}{4}} + T^\varepsilon N^{\frac{7}{12}}K^{\frac{5}{4}} + T^{\frac{1}{4}+\varepsilon} N^{\frac{7}{8}} K^{\frac{1}{4}} \\ \ \ \ \ + T^{\frac{1}{2}+\varepsilon}N^{\frac{5}{12}}K^{\frac{1}{4}}(N^{\frac{5}{24}} + K^{\frac{1}{4}}),  &\text{ if } \ \ \ a_n =\mu^2(n)(\mu * \Lambda^{*k})f(n) \text{ and } b_n \ll n^\varepsilon, \\
T^{\varepsilon}N^{\frac{1}{4}}K^{\frac{1}{2}}(N^{\frac{3}{4}} + N^{\frac{1}{2}}K^{\frac{1}{2}} + K)+T^{\frac{1}{4}+\varepsilon}N^{\frac{7}{8}}K^{\frac{1}{4}} \\ \ \ \ \ +T^{\frac{1}{2}+\varepsilon}N^{\frac{3}{8}}K^{\frac{1}{4}}(N+K)^{\frac{1}{4}},  &\text{ if } \ \ \ a_n =\mu(n)f(n) \text{ and } b_n \ll n^\varepsilon, \nonumber
\end{cases}
\end{equation}
with $f \in \mathcal{F}$. 
\end{theorem}
\begin{remark}
Theorem \ref{theorem2} states that if one couples the Conrey mollifier along with a generic mollifier under the same twisted second moment, then the limiting exponent in the error $\mathcal{E}$ above becomes $\frac{4}{7}$.
\end{remark}
Lastly, we consider
\[
J(\alpha,\beta) := \int_{-\infty}^\infty \zeta(\tfrac{1}{2}+\alpha+it)\zeta(\tfrac{1}{2}+\beta-it) |A(\tfrac{1}{2}+it)|^2 |B(\tfrac{1}{2}+it)|^2 \Phi \bigg(\frac{t}{T}\bigg)dt. 
\]
where $A$ and $B$ are defined in \eqref{2dirichletseries} and the additional condition that $N \ge K$. The result is as follows.

\begin{theorem} \label{theorem3}
Let $\alpha, \beta \ll L^{-1}$ and $a_n, b_n \ll n^\varepsilon$. Then one has
\begin{align}
  I(\alpha ,\beta ) = \sum\limits_{d,e \le NK}  \frac{\mathfrak{a}_d \overline{\mathfrak{a}}_e}{[d,e]}\frac{(d,e)^{\alpha  + \beta }}{d^\alpha e^\beta }\int_{ - \infty }^\infty   \bigg( \zeta (1 + \alpha  + \beta ) + \zeta (1 - \alpha  - \beta )\bigg( \frac{2\pi de}{t(d,e)^2} \bigg)^{\alpha  + \beta } \bigg)\Phi\bigg(\frac{t}{T}\bigg)dt  + O(\mathcal{E}), \nonumber
\end{align}
where $\mathcal{E}$ given by
\begin{align}
\mathcal{E} = T^\varepsilon(T^{\frac{1}{2}}N^{\frac{3}{4}}K+T^{\frac{1}{2}}NK^{\frac{1}{2}}+N^{\frac{7}{4}}K^{\frac{3}{2}}), \nonumber
\end{align}
and with $\mathfrak{a}_d : =\sum_{nk=d} a_n b_k$.
\end{theorem}

\subsection{Final remarks}
The original approach presented in \cite{conrey89} to get the main terms of $I$ required the \textsl{functional equation} of the more complicated Estermann zeta-function (at $u=0$) \cite{estermann},
\[
E_u(s, a/q) := \sum_{n=1}^\infty \frac{\sigma_u(n)\e (an/q)}{n^s}, \quad \textnormal{for} \quad \real(s) > \real(u)+1, \quad a,q \in \N \textnormal{ such that } (a,q)=1.
\]
The more modern method of \cite{bcr, bcy, youngsimple} utilizes the \textsl{approximate functional equation} of the simpler Riemann zeta-function. Both techniques are equivalent, in the sense that they lead to the same main and error terms. We have shown preference for the latter in order to parallel \cite{bcr}.\\

Finally, throughout the paper, we use shall use the convention that $\varepsilon$ denotes and arbitrarily small positive quantity that may not be the same at each occurrence.
\section{Preliminary results}
\subsection{The approximate functional equation}
 As mentioned in the introduction, the starting point is an adaptation of the approximate functional equation of the Riemann zeta function (see \cite[Lemma 4.1]{bcy} and \cite[Lemma 4]{youngsimple}, or more generally \cite[Theorem 5.3]{iwaniecKowalski}). More precisely, let
\[
G(s) := e^{s^2}p(s) \quad \textnormal{where} \quad p(s) := \frac{(\alpha+\beta)^2 - (2s)^2}{(\alpha+\beta)^2}.
\]
In other words, $G$ is an entire function such that $G(x+iy) \ll_x y^{-A}$ for any fixed $x$ and $A>0$. We note that $G(-s)=G(s)$, $G(0)=1$ and $G(\pm \frac{\alpha+\beta}{2})=0$. Next, we define
\begin{align} \label{defV}
W(x) := \frac{1}{2\pi i} \int_{(2)}x^{-w}G(w)\frac{dw}{w} \quad \textnormal{and} \quad V_{\alpha ,\beta }(x,t) := \frac{1}{2\pi i}\int_{(2)}  \frac{G(s)}{s}{g_{\alpha ,\beta }}(s,t)x^{ - s}ds
\end{align}
where
\begin{align} \label{gasymp}
g_{\alpha ,\beta}(s,t) := \pi^{ - s}\frac{\Gamma (\tfrac{1/2 + \alpha  + s + it}{2})\Gamma (\tfrac{1/2 + \beta  + s - it}{2})}{\Gamma (\tfrac{1/2 + \alpha  + it}{2})\Gamma (\tfrac{1/2 + \beta  - it}{2})} = {\bigg( \frac{t}{2\pi} \bigg)^s}(1 + O(t^{ - 1}(1 + |s|^2))),
\end{align}
for large $t$ and $s$ in any fixed vertical strip. With this notation, the approximate functional equation becomes
\begin{align} \label{AFE}
  \zeta \bigg( {\frac{1}{2} + \alpha  + it} \bigg)\zeta \bigg( \frac{1}{2} + \beta  - it \bigg) &= \sum\limits_{m_1,m_2}  \frac{1}{m_1^{1/2 + \alpha }m_2^{1/2 + \beta }}\bigg( \frac{m_1}{m_2} \bigg)^{it}V_{\alpha ,\beta }(m_1 m_2,t) \nonumber \\
   &\quad + X_{\alpha ,\beta ,t}\sum\limits_{m_1, m_2}  \frac{1}{m_1^{1/2 - \beta }m_2^{1/2 - \alpha }}\bigg( \frac{m_1}{m_2} \bigg)^{ it}V_{ - \beta , - \alpha }(m_1 m_2,t) \nonumber \\
	 &\quad + O_A((1 + |t|)^{ - A}),
\end{align}
for $\alpha, \beta$ with real part less than $1/2$, and for any $A \ge 0$. Here
\begin{align} \label{Xdefandasymp}
X_{\alpha ,\beta ,t} := \pi^{\alpha  + \beta }\frac{\Gamma (\tfrac{1/2 - \alpha  - it}{2})\Gamma (\tfrac{1/2 - \beta  + it}{2})}{\Gamma (\tfrac{1/2 + \alpha  + it}{2})\Gamma (\tfrac{1/2 + \beta  - it}{2})} = \bigg( \frac{t}{2\pi } \bigg)^{ - \alpha  - \beta }(1 + O(t^{ - 1})),
\end{align}
for large $t$ and $s$ in any fixed vertical strip. As remarked in \cite[p. 43]{bcy}, $G(s)$, also known as a pole annihilator, can be chosen from a wide class of functions. This particular choice has the advantage of making $G$ vanish at $s = \pm \frac{\alpha+\beta}{2}$. If $\alpha, \beta \to 0$, then \eqref{AFE} becomes
\begin{align}
  \left| \zeta \bigg( {\frac{1}{2} + it} \bigg) \right|^2 = 2\sum\limits_{m_1,m_2}  \frac{1}{(m_1m_2)^{1/2}}\bigg( \frac{m_1}{m_2} \bigg)^{it}V_{0,0}(m_1m_2,t) + O_A((1 + |t|^{ - A}) \nonumber 
\end{align}
and we also get
\begin{align}
  V_{0,0}({m_1}{m_2},t) &= \frac{1}{{2\pi i}}\int_{(2)} {} \frac{{G(s)}}{s}{g_{0,0}}(s,t){({m_1}{m_2})^{ - s}}ds = \frac{1}{{2\pi i}}\int_{(2)} {} \frac{{G(s)}}{s}{\bigg( {\frac{{2\pi {m_1}{m_2}}}{t}} \bigg)^{ - s}}ds \nonumber \\
	 &\quad + O\bigg( \frac{1}{t^{1/2-\varepsilon}(m_1 m_2)^{1/2+\varepsilon}}\bigg) \nonumber \\
   &= W\bigg( {\frac{{2\pi {m_1}{m_2}}}{t}} \bigg) + O\bigg( \frac{1}{t^{1/2-\varepsilon}(m_1 m_2)^{1/2+\varepsilon}}\bigg). \nonumber  
\end{align}
Therefore
\[
\left| {\zeta \bigg( {\frac{1}{2} + it} \bigg)} \right|^2 = 2\sum\limits_{{m_1},{m_2}} {} \frac{1}{{{{({m_1}{m_2})}^{1/2}}}}{\bigg( {\frac{{{m_1}}}{{{m_2}}}} \bigg)^{it}}W\bigg( {\frac{{2\pi {m_1}{m_2}}}{t}} \bigg) + O(T^{-1/2+\varepsilon}).
\]
This is, in fact, the starting point of \cite{bcr}.
\subsection{Bounds on Kloosterman sums}
We require some results on Kloosterman sums that will be used to bound various error terms. We start with a result of Deshouillers and Iwaniec \cite{deshouillersIwaniec1,deshouillersIwaniec2}.
\begin{lemma}\label{DIlemma}
Assume $A,B,N,V \geq 1$ and $|c(a,n)| \leq 1$. Then
\begin{align*}
&\sumtwo_{\substack{v \leq V, b \leq B \\ (b\varrho,v)=1}} \bigg| \sum_{n \leq N} \sum_{\substack{a \leq A \\ (a,v)=1}} c(a,n) \e \bigg(n\frac{ \overline{\varrho ab}}{v}\bigg) \bigg| \\
&\ll (ABNV)^{1/2+\epsilon} \{(BV)^{1/2} + (A+N)^{1/4} [BV(N + \varrho A)(V + \varrho A^2) + \varrho A^2B^2 N]^{1/4} \}.
\end{align*}
\end{lemma}

The next result is due to Bettin and Chandee \cite{bc}, and it improves a result of Duke, Friedlander and Iwaniec \cite{dfi} on bounds of bilinear Kloosterman sums.
\begin{lemma}\label{BClemma}
Let $\alpha_m, \beta_n, \nu_a$ be complex numbers, where $M \le m < 2M$, $N \le n < 2N$, and $A \le a < 2A$. Then for any $\varepsilon > 0$, we have
\begin{align*}
\sum_{a \sim A} \sumtwo_{\substack{m \sim M, n \sim N \\ (m,n)=1}} \nu_a \alpha_m \beta_n \e \bigg(\frac{a\overline{m}}{n}\bigg) &\ll_\varepsilon \| \alpha\| \|\beta\| \|\nu\| \bigg(1+\frac{A}{MN}\bigg)^{\frac{1}{2}} \\
& \quad \times ((AMN)^{\frac{7}{20}+\varepsilon}(M+N)^{\frac{1}{4}}+(AMN)^{\frac{3}{8}+\varepsilon}(AN+AM)^{\frac{1}{8}}),
\end{align*}
where $\|\cdot\|$ denotes the $L_2$ norm.
\end{lemma}

We may now proceed with the proof of the results.
\section{Proof of Theorem \ref{theorem1}}
From \eqref{AFE} we get
\[
I(\alpha ,\beta ) = {I_1}(\alpha ,\beta ) + {I_2}(\alpha ,\beta ) + O(T^{-A}),
\]
where
\begin{align}
I_1(\alpha ,\beta ) &= \int_{ - \infty }^\infty  {} \sum\limits_{{m_1},{m_2}} {} \frac{1}{{m_1^{1/2 + \alpha }m_2^{1/2 + \beta }}}{\bigg( {\frac{{{m_1}}}{{{m_2}}}} \bigg)^{  it}}{V_{\alpha ,\beta }}({m_1}{m_2},t)\sum\limits_{{n_1} \le N} {} \frac{{{a_{{n_1}}}}}{{n_1^{1/2 + it}}}\sum\limits_{{n_2} \le N} {} \frac{{{{\overline{a}}_{{n_2}}}}}{{n_2^{1/2 - it}}}\Phi \bigg( {\frac{t}{T}} \bigg)dt \nonumber \\
&= \sum_{\substack{m_1,m_2,n_1,n_2 }} \frac{{{a_{{n_1}}}{{\bar a}_{{n_2}}}}}{{m_1^{1/2 + \alpha }m_2^{1/2 +\beta }n_1^{1/2}n_2^{1/2}}}\int_{ - \infty }^\infty   {\left( {\frac{{{m_1}{n_2}}}{{{m_2}{n_1}}}} \right)^{it}}{V_{\alpha, \beta }}({m_1}{m_2},t)\Phi \left( {\frac{t}{T}} \right)dt, \nonumber
\end{align}
and
\begin{align}
  I_2(\alpha ,\beta ) 
   &= \sum_{\substack{m_1,m_2,n_1,n_2 }} \frac{{{a_{{n_1}}}{{\bar a}_{{n_2}}}}}{{m_1^{1/2 - \beta }m_2^{1/2 - \alpha }n_1^{1/2}n_2^{1/2}}}\int_{ - \infty }^\infty   {\left( {\frac{{{m_1}{n_2}}}{{{m_2}{n_1}}}} \right)^{it}}{V_{ - \beta , - \alpha }}({m_1}{m_2},t)X_{\alpha,\beta,t}\Phi \left( {\frac{t}{T}} \right)dt. \nonumber  
\end{align}

We first concentrate on $I_1$, then describe the modifications necessary to handle $I_2$. Pulling the sums out of the integrals, we get
\begin{align}
  {I_1}(\alpha ,\beta ) &= \sum\limits_{{m_1},{m_2},{n_1},{n_2}} {} \frac{{{a_{{n_1}}}{{\overline{a}}_{{n_2}}}}}{{m_1^{1/2 + \alpha }m_2^{1/2 + \beta }n_1^{1/2}n_2^{1/2}}}\int_{ - \infty }^\infty  {} {\bigg( {\frac{{{m_1}{n_2}}}{{{m_2}{n_1}}}} \bigg)^{  it}}{V_{\alpha ,\beta }}({m_1}{m_2},t)\Phi \bigg( {\frac{t}{T}} \bigg)dt \nonumber \\
   &= \mathcal{D}_1 + \mathcal{S}_1, \nonumber  
\end{align}
where the sum is over $n_1, n_2 \le N$, $\mathcal{D}_1$ is the sum when $m_1n_2 = m_2n_1$ and $\mathcal{S}_1$ is the sum when $m_1n_2 \ne m_2n_1$.

\subsection{Diagonal terms} We start with $\mathcal{D}_1$. For $j=1,2$, we write $m_j = \ell n_j^*$ where $n_j^* = \frac{n_j}{(n_1,n_2)}$. The contribution from the diagonal term is therefore
\begin{align}
  \mathcal{D}_1 &= \sum\limits_{{m_1},{m_2},{n_1},{n_2}} {} \frac{{{a_{{n_1}}}{{\overline{a}}_{{n_2}}}}}{{m_1^{1/2 + \alpha }m_2^{1/2 + \beta }n_1^{1/2}n_2^{1/2}}}\int_{ - \infty }^\infty  {} {V_{\alpha ,\beta }}({m_1}{m_2},t)\Phi \bigg( {\frac{t}{T}} \bigg)dt \nonumber \\
   &= \sum\limits_{\ell ,{n_1},{n_2}} {} \frac{{{a_{{n_1}}}{{\overline{a}}_{{n_2}}}{{({n_1},{n_2})}^{1 + \alpha  + \beta }}}}{{{\ell ^{1 + \alpha  + \beta }}n_1^{1 + \alpha }n_2^{1 + \beta }}}\int_{ - \infty }^\infty  {} {V_{\alpha ,\beta }}(\ell n_1^*\ell n_2^*,t)\Phi \bigg( {\frac{t}{T}} \bigg)dt \nonumber \\
   &= \sum\limits_{\ell ,{n_1},{n_2}} {} \frac{{{a_{{n_1}}}{{\overline{a}}_{{n_2}}}{{({n_1},{n_2})}^{1 + \alpha  + \beta }}}}{{{\ell ^{1 + \alpha  + \beta }}n_1^{1 + \alpha }n_2^{1 + \beta }}}\int_{ - \infty }^\infty  {} \frac{1}{{2\pi i}}\int_{(2)} {} {g_{\alpha ,\beta }}(s,t){({\ell ^2}n_1^*n_2^*)^{ - s}}G(s)\frac{{ds}}{s}\Phi \bigg( {\frac{t}{T}} \bigg)dt \nonumber \\
   &= \sum\limits_{{n_1},{n_2}} {} \frac{{{a_{{n_1}}}{{\overline{a}}_{{n_2}}}{{({n_1},{n_2})}^{1 + \alpha  + \beta }}}}{{n_1^{1 + \alpha }n_2^{1 + \beta }}}\int_{ - \infty }^\infty  {} \frac{1}{{2\pi i}}\int_{(2)} {} {g_{\alpha ,\beta }}(s,t){(n_1^*n_2^*)^{ - s}}\zeta (1 + \alpha  + \beta  + 2s)G(s)\frac{{ds}}{s}\Phi \bigg( {\frac{t}{T}} \bigg)dt \nonumber \\
   &= \sum\limits_{{n_1},{n_2}} {} \frac{{{a_{{n_1}}}{{\overline{a}}_{{n_2}}}{{({n_1},{n_2})}^{1 + \alpha  + \beta }}}}{{n_1^{1 + \alpha }n_2^{1 + \beta }}}\int_{ - \infty }^\infty  {} \frac{1}{{2\pi i}}\int_{(2)} {} {\bigg( {\frac{{2\pi n_1^*n_2^*}}{t}} \bigg)^{ - s}}\zeta (1 + \alpha  + \beta  + 2s)G(s)\frac{{ds}}{s}\Phi \bigg( {\frac{t}{T}} \bigg)dt \nonumber \\
	 &\quad + O(T^{1/2+\varepsilon}), \nonumber  
\end{align}
by the use of \eqref{gasymp} in the last line. This term will be later combined with a contribution from the off-diagonal terms. Together, they give the main term in \ref{theorem1}. 
\subsection{Off-Diagonal terms}
Let us now move on to the off-diagonal terms $\mathcal{S}_1$. First, recall that
\begin{align}
  \mathcal{S}_1 = \sum_{\substack{m_1,m_2,n_1,n_2 \\ m_1 n_2 \ne m_2n_1}}  \frac{{{a_{{n_1}}}{{\overline{a}}_{{n_2}}}}}{{m_1^{1/2 + \alpha }m_2^{1/2 + \beta }n_1^{1/2}n_2^{1/2}}}\int_{ - \infty }^\infty  {} {\bigg( {\frac{{{m_1}{n_2}}}{{{m_2}{n_1}}}} \bigg)^{it}}{V_{\alpha ,\beta }}({m_1}{m_2},t)\Phi \bigg( {\frac{t}{T}} \bigg)dt .  \nonumber 
\end{align} 
We now write $m_1 n_2 - m_2 n_1 = \Delta$. From (4.4) of \cite{youngsimple} (note the typo) or \cite[Proposition 5.4]{iwaniecKowalski} we have that for any $A \ge 0$ and $j=0,1,2,\cdots$, we have uniformly in $x$,
\begin{align} \label{derivativeboundV}
t^j \frac{\partial^j}{\partial t^j} V_{\alpha,\beta}(x,t) \ll_{A,j} (1+|x/t|)^{-A}.
\end{align}
This means that we can truncate the sum over $m_1, m_2$ to $m_1 m_2 \le T^{1+\varepsilon}$. We introduce the smooth partition of unity
\[
1 = \sideset{}{'}\sum_{M} F_M(x), \quad T^{-100} \le x \le T^{1+\varepsilon},
\]
where $F_M(x)$ is smooth, supported in $[M/2,3M]$, and satisfies $F_M^{(j)}(x) \ll_j M^{-j}$ for all $j \ge 0$. This partition of unity will also satisfy $\sideset{}{'}\sum\nolimits_{M} 1 \ll \log (2+T)$. Therefore the non-diagonal term becomes
\begin{align} \label{newmathcalS}
  \mathcal{S}_1 &= \sideset{}{'}\sum\limits_{{N_1}}  \sideset{}{'}\sum\limits_{{N_2}} {} \sideset{}{'}\sum\limits_M  \sum\limits_{\Delta  \ne 0} \sum_{\substack{{m_1},{m_2},{n_1},{n_2} \\ m_1n_2 - m_2n_1 = \Delta}} {\frac{{{a_{{n_1}}}{{\overline{a}}_{{n_2}}}}}{{m_1^{1/2 + \alpha }m_2^{1/2 + \beta }n_1^{1/2}n_2^{1/2}}}}  \nonumber \\
   &\quad \times \bigg( {\int_{ - \infty }^\infty  {} {{\bigg( {1 + \frac{\Delta }{{{m_2}{n_1}}}} \bigg)}^{it}} V_{\alpha,\beta}(m_1m_2,t)\Phi \bigg( {\frac{t}{T}} \bigg)dt} \bigg){F_{{N_1}}}({n_1}){F_{{N_2}}}({n_2}){F_M}({m_2}) + O(1),
\end{align}
where $N_1, N_2 \le N$ and $M \le T^{1+\varepsilon}$. Using \eqref{derivativeboundV} and $\ell$ integration by parts, it is not hard to show that with $|\Delta| > D$, where $D:=\frac{MN_1}{T^{1-\varepsilon}}$, give a negligible $\ll_{A,\varepsilon} T^{-A}$ contribution. In other words
\begin{align}
\sideset{}{'}\sum_{N_1, N_2,M} &\sum_{|\Delta|>D} \sum_{n_1, m_2}\sum_{\substack{n_2,m_1 \\ m_1n_2 - m_2n_1 = \Delta}} \frac{a_{n_1}\bar{a}_{n_2}}{m_1^{1/2+\alpha}m_2^{1/2+\beta}n_1^{1/2}n_2^{1/2}} \nonumber \\
& \quad \times \int_{-\infty}^\infty \bigg(1+\frac{\Delta}{m_2n_1}\bigg)^{it} V_{\alpha,\beta}(m_1m_2,t) \Phi\bigg(\frac{t}{T}\bigg)dt \nonumber \\
& \ll_{\ell} \sideset{}{'}\sum_{N_1,M} \sum_{D < |\Delta| \leq T^{O(1)}} \sum_{\substack{n_1 \sim N_1 \\ m_2 \sim M}} \sum_{\substack{m_1 \le T^{1+\varepsilon} \\ n_2 \le N \\ m_1 n_2 - m_2 n_1 = \Delta}} \frac{T^{-\ell+1+\varepsilon}}{m_1^{1/2}m_2^{1/2}n_1^{1/2}n_2^{1/2}} \bigg|\log\bigg(1+\frac{\Delta}{m_2n_1}\bigg)\bigg|^{-\ell} \nonumber \\
& \ll_{\ell} \sideset{}{'}\sum_{N_1,M} \sum_{D < |\Delta| \leq T^{O(1)}} \sum_{\substack{n_1 \sim N_1 \\ m_2 \sim M}} \frac{1}{\sqrt{m_2n_1}} \sum_{\substack{m_1 \le T^{1+\varepsilon} \\ n_2 \le N \\ m_1 n_2 - m_2 n_1 = \Delta}} \frac{1}{\sqrt{m_1n_2}} T^{-\ell+1+\varepsilon} \bigg(\frac{MN_1}{D}\bigg)^\ell \ll_{A,\varepsilon} T^{-A}, \nonumber
\end{align}
where $\ell$ is large enough. Next, we move on to $|\Delta| < D$. By \eqref{gasymp} we have
\[
V_{\alpha,\beta}(x,t) = W\bigg(\frac{2\pi x}{t}\bigg) + \frac{1}{2 \pi i} \int_{(2)} \frac{G(s)}{s} E(s,t)x^{-s}ds,
\]
where $E(s,t)$ is an analytic function of $s$ and $t$ for $t$ sufficiently large and $\real(s)>0$. Moreover, we have the estimate
\[
E(s,t) \ll_\sigma \frac{1+|s|^2}{t^{1-\sigma}}.
\]
The error term associated to this approximation is then given by
\begin{align*}
\mathcal{E}_E &:= \sideset{}{'}\sum\limits_{{N_1}}  \sideset{}{'}\sum\limits_{{N_2}} {} \sideset{}{'}\sum\limits_M  \sum\limits_{0 < |\Delta| \leq D} \sum_{\substack{{m_1},{m_2},{n_1},{n_2} \\ {m_1}{n_2} - {m_2}{n_1} = \Delta}} \frac{{{a_{{n_1}}}{{\overline{a}}_{{n_2}}}}}{{m_1^{1/2 + \alpha }m_2^{1/2 + \beta }n_1^{1/2}n_2^{1/2}}} {F_{{N_1}}}({n_1}){F_{{N_2}}}({n_2}){F_M}({m_2})  \nonumber \\
   &\quad \times \bigg( {\int_{ - \infty }^\infty  {} {{\bigg( {1 + \frac{\Delta }{{{m_2}{n_1}}}} \bigg)}^{it}}\bigg( \frac{1}{2\pi i} \int_{(2)} \frac{G(s)}{s} E(s,t) (m_1m_2)^{-s}ds \bigg)\Phi \bigg( {\frac{t}{T}} \bigg)dt} \bigg).
\end{align*}
As $E(s,t)$ is analytic in $s$ for $\sigma > 0$, we can move the line of integration in $s$ from $\sigma = 2$ to $\sigma = \varepsilon$. We now apply the triangle inequality and proceed to use trivial estimations. We upper bound the quantity $(m_1m_2)^{-\varepsilon}$ by 1, and use our bound for $E(s,t)$. From the rapid decay of $G(s)$ in vertical strips and the fact that $\Phi$ is supported in $[1,2]$, it is then easy to see that
\begin{align*}
&\bigg|\int_{ - \infty }^\infty   \bigg( {1 + \frac{\Delta }{{{m_2}{n_1}}}} \bigg)^{it}\bigg( \frac{1}{2\pi i} \int_{(\varepsilon)} \frac{G(s)}{s} E(s,t) (m_1m_2)^{-s}ds \bigg)\Phi \bigg( {\frac{t}{T}} \bigg)dt\bigg| \\
&\ll_\varepsilon \frac{T^\varepsilon}{T} \int_{-\infty}^\infty \Phi \left( \frac{t}{T} \right) \int_{(\varepsilon)} \frac{|G(s)|}{|s|} (1+|s|^2) |ds| \ll_\varepsilon T^\varepsilon.
\end{align*}
We deduce that
\begin{align*}
\mathcal{E}_E &\ll_\varepsilon T^\varepsilon \sideset{}{'}\sum\limits_{{N_1}}  \sideset{}{'}\sum\limits_{{N_2}} {} \sideset{}{'}\sum\limits_M  \sum\limits_{0 < |\Delta| \leq D} \sum_{\substack{{m_1},{m_2},{n_1},{n_2} \\ {m_1}{n_2} - {m_2}{n_1} = \Delta}} {\frac{{F_{{N_1}}}({n_1}){F_{{N_2}}}({n_2}){F_M}({m_2})}{{m_1^{1/2 }m_2^{1/2 }n_1^{1/2}n_2^{1/2}}}}.
\end{align*}
Since $|\Delta| \ll D$ we have $\frac{\Delta}{m_2n_1} \ll \frac{1}{T^{1-\varepsilon}}$, and
\begin{align*}
\frac{1}{m_1} &= \frac{n_2}{m_2n_1} \left(1 + \frac{\Delta}{m_2n_1} \right)^{-1} \ll \frac{n_2}{m_2n_1}.
\end{align*}
Taking square roots and using the same bounds and approximations as in \cite{bcr} (see also the treatment below), we then obtain
\begin{align*}
\mathcal{E}_E &\ll T^\varepsilon \sideset{}{'}\sum\limits_{{N_1}}  \sideset{}{'}\sum\limits_{{N_2}} {} \sideset{}{'}\sum_{M \ll T^{1/2 + \varepsilon} \sqrt{\frac{N_2}{N_1}}}  \sum\limits_{0 < |\Delta| \leq D} \sum_{\substack{{m_2},{n_1},{n_2} \\ m_2n_1 \equiv - \Delta \pmod{n_2}}} {\frac{{F_{{N_1}}}({n_1}){F_{{N_2}}}({n_2}){F_M}({m_2})}{m_2n_1}      }.
\end{align*}
Next, we extract the greatest common divisor $d$ of $n_1,n_2$, and obtain
\begin{align*}
\mathcal{E}_E &\ll T^\varepsilon \sum_{d\leq N} \frac{1}{d} \sideset{}{'}\sum\limits_{{N_1}}  \sideset{}{'}\sum\limits_{{N_2}} {} \sideset{}{'}\sum_{M \ll T^{1/2 + \varepsilon} \sqrt{\frac{N_2}{N_1}}} \sum_{0 < |\Delta| \leq \frac{D}{d}}\sum_{(n_1,n_2)=1} \frac{F_{N_1}(dn_1) F_{N_2}(dn_2)}{n_1} \sum_{m_2 \equiv -\overline{n_1}\Delta (n_2)} \frac{F_M(m_2)}{m_2}.
\end{align*}
We study the innermost sum. By the support of $F_M$, the inner sum is bounded by
\begin{align*}
\sum_{\substack{c_1 M < m_2 \leq c_2 M \\ m_2 \equiv v (n_2)}} \frac{1}{m_2},
\end{align*}
for some positive constants $c_1 < c_2$ and some residue class $v$ modulo $n_2$. We change variables and approximate the sum by an integral, which yields
\begin{align*}
\sum_{\substack{c_1 M < m_2 \leq c_2 M \\ m_2 \equiv v \modu (n_2)}} \frac{1}{m_2} \ll \frac{1}{M} + \frac{1}{N_2}.
\end{align*}
Summing over the rest of the variables trivially, we obtain
\begin{align*}
\mathcal{E}_E \ll \frac{T^\varepsilon}{T} M N_1 N_2 \left( \frac{1}{M} + \frac{1}{N_2} \right),
\end{align*}
for some $M \ll T^{1/2 + \varepsilon} \sqrt{\frac{N_2}{N_1}}$ and $N_1,N_2 \leq N$. Applying these bounds yields
\begin{align*}
\mathcal{E}_E &\ll T^\varepsilon \left(\frac{N}{T^{1/2}} + \frac{N^2}{T} \right) \ll N^{1/2} T^{1/2 + \varepsilon}.
\end{align*}

In the case $|\Delta| < D$ we have $\frac{\Delta}{m_2 n_1} \ll \frac{1}{T^{1-\varepsilon}}$, and
\[
m_1 = \frac{m_2 n_1 + \Delta}{n_2} = m_2 \frac{n_1}{n_2} \bigg( 1 + \frac{\Delta}{m_2 n_1} \bigg).
\]
Hence for $T < t < 2T$, we get the approximation
\begin{align}
\bigg(1 + \frac{\Delta}{mn}\bigg)^{it} &= e^{it \log (1 + \tfrac{\Delta}{m_2 n_1})} = e^{it \tfrac{\Delta}{m_2 n_1}} \bigg( 1 - \frac{it\Delta^2}{2 m_2^2 n_1^2} + O \bigg( \frac{1}{T^{2-\varepsilon}}\bigg)\bigg), \nonumber
\end{align}
as well as 
\[
W\bigg(\frac{2\pi m_1 m_2}{t}\bigg) = W\bigg(\frac{2\pi m_2^2 n_1}{tn_2}\bigg) + \frac{2 \pi m_2 \Delta}{t n_2} W'\bigg(\frac{2\pi m_1 m_2}{t}\bigg) + O \bigg(\frac{1}{T^{2-\varepsilon}}\bigg).
\]
Since $m_1 m_2 \le T^{1+\varepsilon}$, we have that $m_2(m_2 n_1 + \Delta) \le n_2 T^{1+\varepsilon}$. So, we have the bound $M \ll T^{1/2+\varepsilon}\sqrt{\frac{N_2}{N_1}}$, and the error term from using the above approximations in \eqref{newmathcalS} is
\begin{align}
&\ll \frac{T}{T^{2-\varepsilon}} \sideset{}{'}\sum\limits_{N_1}\sideset{}{'}\sum\limits_{N_2} \sideset{}{'}\sum_{M \ll T^{1/2+\varepsilon}\sqrt{N_2}{N_1}} \sum_{0 \le |\Delta| \le D} \sum_{\substack{n_1 \sim N_1 \\ n_2 \sim N_2}} \sum_{m_2 \sim M} \frac{1}{m_2n_1} \nonumber \\
&\ll \sideset{}{'}\sum\limits_{N_1}\sideset{}{'}\sum\limits_{N_2} \sideset{}{'}\sum_{M \ll T^{1/2+\varepsilon}\sqrt{N_2}{N_1}} \frac{MN_1N_2}{T^{2-\varepsilon}} \ll \frac{\sqrt{T}N_2^{3/2}N_1^{1/2}}{T^{2-\varepsilon}} \ll \frac{N^2}{T^{3/2-\varepsilon}}, \nonumber
\end{align}
by the definition of $D$ and the above bound on $M$. This gives us the split
\begin{align} \label{splitSAE}
\mathcal{S}_1 = \mathcal{A}_1 + \mathcal{E}_W + O \bigg( 1 + \frac{N^2}{T^{3/2-\varepsilon}}\bigg),
\end{align}
where
\begin{align} \label{defmathcalA}
\mathcal{A}_1 &= \sideset{}{'}\sum\limits_{N_1} \sideset{}{'}\sum\limits_{N_2}  \sideset{}{'}\sum\limits_{M \le T^{1/2+\varepsilon} \sqrt{\frac{N_2}{N_1}}} \sum_{0 < |\Delta| \le D} \sum_{n_1, n_2} \sum_{\substack{n_1 m_2 \equiv -\Delta {\kern 1pt} (\modu n_2) \\ m_2 > 0}} \frac{a_{n_1}\bar{a}_{n_2}}{m_2^{1+\alpha+\beta} n_1} \bigg( \frac{n_2}{n_1} \bigg)^\alpha \nonumber \\
& \quad \times \bigg( \int_{-\infty}^\infty \e\bigg( \frac{\Delta t}{2 \pi m_2 n_1}\bigg) W\bigg(\frac{2\pi m_2^2 n_1}{tn_2}\bigg) \Phi \bigg( \frac{t}{T}\bigg) dt\bigg) F_M(m_2)F_{N_1}(n_1)F_{N_2}(n_2),
\end{align}
as well as
\begin{align} \label{definitionofE}
\mathcal{E}_W & = \sideset{}{'}\sum\limits_{N_1,N_2} \sideset{}{'}\sum\limits_{M \le T^{1/2+\varepsilon} \sqrt{\frac{N_2}{N_1}}} \sum_{0 < |\Delta| \le D} \sum_{n_1, n_2} \sum_{\substack{n_1 m_2 \equiv -\Delta {\kern 1pt} (\modu n_2) \\ m_2 > 0}} \frac{a_{n_1}\bar{a}_{n_2}}{m_2^{1+\alpha+\beta} n_1} \bigg( \frac{n_2}{n_1} \bigg)^\alpha F_M(m_2)F_{N_1}(n_1)F_{N_2}(n_2) \nonumber \\
& \quad \times  \int_{-\infty}^\infty \e\bigg( \frac{\Delta t}{2 \pi m_2 n_1}\bigg) \bigg[ W\bigg(\frac{2\pi m_2^2 n_1}{tn_2}\bigg) \bigg( -\frac{(\tfrac{1}{2}+\alpha)\Delta}{m_2 n_1}-\frac{it\Delta^2}{2 m_2^2 n_1^2} \bigg) +  \frac{2 \pi m_2 \Delta}{t n_2} W'\bigg(\frac{2\pi m_2^2 n_1}{tn_2}\bigg) \bigg] d\Phi ,
\end{align}
since the rest of the terms arising from the above approximations also give a contribution which is $\ll N^2 T^{-3/2+\varepsilon}$. We set the temporary notation $d\Phi := \Phi(\frac{t}{T})dt$. Let us first examine $\mathcal{A}_1$. To this end, extract the common divisor $d$ of $n_1$ and $n_2$ and write $\mathcal{A}$ as
\[
\mathcal{A}_1 = \sum_{d \le N} \frac{1}{d} \sideset{}{'}\sum\limits_{N_1, N_2 \le N} \sideset{}{'}\sum\limits_{M \le T^{1/2+\varepsilon} \sqrt{\frac{N_2}{N_1}}} \sum_{0<|\Delta|\le \frac{D}{d}}\sum_{\substack{n_1, n_2 \\ (n_1,n_2)=1}} a_{dn_1} \bar{a}_{dn_2} F_{N_1}(dn_1) F_{N_2}(dn_2) \mathcal{A}_{M,N_i}(n_1, n_2, \Delta),
\]
where
\[
\mathcal{A}_{1,M,N_i}({n_1},{n_2},\Delta ) = \sum\limits_{{m_2} \equiv  - \overline {{n_1}} \Delta {\kern 1pt} (\modu n_2)} {} \frac{{{F_M}({m_2})}}{{m_2^{1 + \alpha  + \beta }{n_1}}}{\bigg( {\frac{{{n_2}}}{{{n_1}}}} \bigg)^\alpha }\int_{ - \infty }^\infty  {\e\bigg( {\frac{{\Delta t}}{{2\pi {m_2}{n_1}}}} \bigg)W\bigg(\frac{2\pi m_2^2 n_1}{tn_2}\bigg)\Phi \bigg( {\frac{t}{T}} \bigg)dt} .
\]
Using Poisson's summation formula yields
\begin{align}
  \mathcal{A}_{1,M,N_i}(n_1,n_2,\Delta ) &= \frac{1}{n_1n_2} \bigg( {\frac{{{n_2}}}{{{n_1}}}} \bigg)^\alpha \sum\limits_{h \in \mathbb{Z}} {\operatorname{e} \bigg( { - \frac{{h\overline {{n_1}} \Delta }}{{{n_2}}}} \bigg)} \int_0^\infty  {\operatorname{e} \bigg( { - \frac{{hx}}{{{n_2}}}} \bigg)\frac{{{F_M}(x)}}{{{x^{1 + \alpha  + \beta }}}}}  \nonumber \\
   &\quad \times \int_{ - \infty }^\infty \operatorname{e} \bigg( \frac{\Delta t}{2\pi xn_1} \bigg)W\left(\frac{2 \pi x^2 n_1}{t n_2}\right)\Phi \bigg( {\frac{t}{T}} \bigg)dt dx. \nonumber 
\end{align}
Next, we make the change of variable $x \to \frac{x}{n_1}$ so that
\[
\mathcal{A}_{1,M,N_i}(n_1,n_2,\Delta ) = \frac{n_1^\beta n_2^\alpha }{n_1n_2} \sum\limits_{h \in \mathbb{Z}}  {{\tilde {\mathcal{A}}}_{M,{N_i}}}(h,{n_1},{n_2},\Delta )\operatorname{e} \bigg( { - \frac{{h\overline {{n_1}} \Delta }}{{{n_2}}}} \bigg),
\]
where
\[
\tilde{\mathcal{A}}_{1,M,N_i}(h,n_1,n_2,\Delta ) = \int_0^\infty  \operatorname{e} \bigg( { - \frac{{hx}}{{{n_1}{n_2}}}} \bigg)\frac{{{F_M}(\tfrac{x}{{{n_1}}})}}{{{x^{1 + \alpha  + \beta }}}} \int_{ - \infty }^\infty \operatorname{e} \bigg( {\frac{{\Delta t}}{{2\pi x}}} \bigg)W\left(\frac{2 \pi x^2 }{n_1 n_2 t}\right)\Phi \bigg( {\frac{t}{T}} \bigg)dtdx .
\]
Now we must study three cases: $h=0$, $|h| \ge H_d:=\frac{N_2}{dM}T^\varepsilon$ and $0< |h| <H_d$. The first case will provide the contribution to the main term, the second case is negligible and the third case will require estimates on Kloosterman sums.
\subsubsection{The case $h=0$.}
The contribution to $\mathcal{A}_1$ from $h=0$ is
\begin{align} \label{defA0}
\mathcal{A}_{1,0} = \sum_{d \le N} &\sideset{}{'}\sumtwo_{N_1 \le N, N_2 \le N} \sideset{}{'}\sum_{M \le T^{1/2+\varepsilon} \sqrt{\frac{N_2}{N_1}}} \sum_{0 < |\Delta| \le \frac{D}{d}} \sum_{\substack{n_1, n_2 \\ (n_1,n_2)=1}} \frac{a_{dn_1}\bar{a}_{dn_2}F_{N_1}(dn_1)F_{N_2}(dn_2)}{d n_1^{1-\beta} n_2^{1-\alpha}} \nonumber \\
& \times \int_{-\infty}^\infty \int_0^\infty F_M \bigg( \frac{x}{n_1} \bigg) \e \bigg( \frac{\Delta t}{2 \pi x} \bigg) W\left(\frac{2 \pi x^2 }{n_1 n_2 t}\right) \frac{dx}{x^{1+\alpha+\beta}} \Phi\bigg( \frac{t}{T}\bigg)dt. 
\end{align}
As argued in \cite{bcr}, we can extend the sum over $\Delta$ to $\Delta \in \mathbb{Z} \backslash \{0\}$, since, as done previously, we can show that the term $|\Delta| \ge D/d$ yields a negligible contribution. Next, we make the change of variables $y=t/x$ and integrate by parts twice the second line of \eqref{defA0} to get
\begin{align}
   & - \frac{1}{{{\Delta ^2}}}\int_{ - \infty }^\infty  {} \int_0^\infty  {} \operatorname{e} \bigg( {\frac{{\Delta y}}{{2\pi }}} \bigg)\frac{{{d^2}}}{{d{y^2}}}\bigg( {{F_M}\bigg( {\frac{t}{{{n_1}y}}} \bigg)W\bigg(\frac{2 \pi t}{n_1n_2 y_2} \bigg)\frac{1}{{{y^{1 - \alpha  - \beta }}}}} \bigg)dy\frac{1}{{{t^{\alpha  + \beta }}}}\Phi \bigg( {\frac{t}{T}} \bigg)dt \nonumber \\
   &=  - \frac{1}{{{\Delta ^2}}}\int_{ - \infty }^\infty  {} \int_{R(t,{n_1})} {} \operatorname{e} \bigg( {\frac{{\Delta y}}{{2\pi }}} \bigg)\frac{{{d^2}}}{{d{y^2}}}\bigg( {{F_M}\bigg( {\frac{t}{{{n_1}y}}} \bigg)W\bigg(\frac{2 \pi t}{n_1n_2 y_2} \bigg)\frac{1}{{{y^{1 - \alpha  - \beta }}}}} \bigg)dy\frac{1}{{{t^{\alpha  + \beta }}}}\Phi \bigg( {\frac{t}{T}} \bigg)dt, \nonumber  
\end{align} 
where $R(t,n_1) = \{ y | T^{-100} < \frac{t}{n_1y} < T^{1/2+\varepsilon}\sqrt{\frac{N_2}{N_1}} \}$ by a trivial estimate of the part of integral over $y$ with $y \in \R_{>0} \backslash R(t,n_1)$ and the properties of $W$ and $F_M$ and $n_1 \ll T$. Now we sum over $M, N_1, N_2$ and $d$. We start with $M$ so that
\begin{align}
& \sideset{}{'} \sum_M \int_{-\infty}^\infty \int_0^\infty F_M \bigg( \frac{x}{n_1} \bigg) \e \bigg( \frac{\Delta t}{2 \pi x} \bigg) W\left(\frac{2 \pi x^2 }{n_1 n_2 t}\right) \frac{dx}{x^{1+\alpha+\beta}} \Phi \bigg( \frac{t}{T} \bigg)dt \nonumber \\
& = -\frac{1}{\Delta^2} \int_{-\infty}^\infty \int_{R(t,n_1)} \e \bigg( \frac{\Delta y}{2 \pi} \bigg) \frac{d^2}{dy^2} \bigg( W\bigg(\frac{2 \pi t}{n_1n_2 y_2} \bigg) \frac{1}{y^{1-\alpha-\beta}}\bigg)dy \frac{1}{t^{\alpha+\beta}} \Phi \bigg( \frac{t}{T} \bigg)dt + O\bigg( \frac{\log(2+T)}{\Delta^2}\bigg) \nonumber \\
& = -\frac{1}{\Delta^2} \int_{-\infty}^\infty \int_0^\infty \e \bigg( \frac{\Delta y}{2 \pi} \bigg) \frac{d^2}{dy^2} \bigg( W\bigg(\frac{2 \pi t}{n_1n_2 y_2} \bigg) \frac{1}{y^{1-\alpha-\beta}}\bigg)dy \frac{1}{t^{\alpha+\beta}} \Phi \bigg( \frac{t}{T} \bigg)dt + O\bigg( \frac{\log(2+T)}{\Delta^2}\bigg). \nonumber  
\end{align}
Summing over $N_1$ and $N_2$ yields
\begin{align}
\mathcal{A}_{1,0} &= -\sum_{d \le N} \sum_{|\Delta| \ne 0} \sum_{\substack{n_1, n_2 \le \frac{N}{d} \\ (n_1, n_2)=1}} \frac{a_{dn_1}\bar{a}_{dn_2}}{d n_1^{1-\beta}n_2^{1-\alpha} \Delta^2} \nonumber \\
              & \quad \times \int_{-\infty}^\infty \int_0^\infty \e \bigg( \frac{\Delta y}{2 \pi} \bigg) \frac{d^2}{dy^2} \bigg( W\bigg(\frac{2 \pi t}{n_1n_2 y_2} \bigg) \frac{1}{y^{1-\alpha-\beta}}\bigg)dy \frac{1}{t^{\alpha+\beta}} \Phi \bigg( \frac{t}{T} \bigg)dt + O(T^\varepsilon) \nonumber \\
              &= -\sum_{d \le N} \sum_{\substack{n_1, n_2 \le \frac{N}{d} \\ (n_1, n_2)=1}} \frac{a_{dn_1}\bar{a}_{dn_2}}{d n_1^{1-\beta}n_2^{1-\alpha} } \nonumber \\
              & \quad \times \int_{-\infty}^\infty \int_0^\infty \sum_{|\Delta| \ne 0} \frac{1}{\Delta^2} \e \bigg( \frac{\Delta y}{2 \pi} \bigg) \frac{d^2}{dy^2} \bigg( W\bigg(\frac{2 \pi t}{n_1n_2 y_2} \bigg) \frac{1}{y^{1-\alpha-\beta}}\bigg)dy \frac{1}{t^{\alpha+\beta}} \Phi \bigg( \frac{t}{T} \bigg)dt + O(T^\varepsilon) \nonumber \\
							& = \mathcal{A}_{1,0,+} + \mathcal{A}_{1,0,-} + O(T^\varepsilon). \nonumber
\end{align}
The interchange of sum over $\Delta$ and the integrals can be justified by absolute convergence. Here $\mathcal{A}_{0,+}$ is the sum over $\Delta>0$ and $\mathcal{A}_{0,-}$ is the sum over $\Delta<0$. Changing variables $\Delta y \to y$ and re-arranging yields
\begin{align}
\mathcal{A}_{1,0,\pm} &= - \sum_{d \le N} \sum_{\substack{n_1, n_2 \le \frac{N}{d} \\ (n_1, n_2)=1}} \frac{a_{dn_1}\bar{a}_{dn_2}}{d n_1^{1-\beta}n_2^{1-\alpha} } \nonumber \\
                    & \quad \times \int_{-\infty}^\infty \int_0^\infty \sum_{\Delta = 1}^\infty \frac{(e^{iy}+e^{-iy})}{\Delta^{\alpha+\beta}} \frac{d^2}{dy^2} \bigg( W\bigg(\frac{2 \pi \Delta^2 t}{n_1n_2 y_2} \bigg) \frac{1}{y^{1-\alpha-\beta}}\bigg)dy \frac{1}{t^{\alpha+\beta}} \Phi \bigg( \frac{t}{T} \bigg)dt \nonumber \\
										&= - \sum_{n_1, n_2 \le N} \frac{a_{n_1}\bar{a}_{n_2}(n_1,n_2)^{1-\alpha-\beta}}{n_1^{1-\beta}n_2^{1-\alpha} } \nonumber \\
										& \quad \times \int_{-\infty}^\infty \int_0^\infty \sum_{\Delta = 1}^\infty \frac{2 \cos y}{\Delta^{\alpha+\beta}}  \frac{d^2}{dy^2} \bigg( W\bigg(\frac{2 \pi t}{n_1^* n_2^* y_2} \bigg) \frac{1}{y^{1-\alpha-\beta}}\bigg)dy \frac{1}{t^{\alpha+\beta}} \Phi \bigg( \frac{t}{T} \bigg)dt, \nonumber 
\end{align}
since we had previously set $n_i^* = \frac{n_i}{(n_1,n_2)}$ for $i=1,2$. The next step is to see that
\begin{align}
  &\sum\limits_{\Delta  = 1}^\infty  {} \frac{1}{{{\Delta ^{\alpha  + \beta }}}}\frac{{{d^2}}}{{d{y^2}}}\bigg( W\bigg(\frac{2 \pi t}{n_1^* n_2^* y_2} \bigg)\frac{1}{{{y^{1 - \alpha  - \beta }}}} \bigg) \nonumber \\
   &= \sum\limits_{\Delta  = 1}^\infty  {} \frac{1}{{{\Delta ^{\alpha  + \beta }}}}\bigg( {\frac{1}{{2\pi i}}\int_{(2)}  G(s) \bigg( \frac{2 \pi \Delta^2 t}{n_1^*n_2^*} \bigg)^{ - s}(2s - 2 + \alpha  + \beta )(2s + \alpha  + \beta  - 1){y^{2s - 3 + \alpha  + \beta }}\frac{{ds}}{s}} \bigg) \nonumber \\
   &= \frac{1}{{2\pi i}}\int_{(\frac{5}{4})} {} \zeta (2s + \alpha  + \beta )\bigg( \frac{2 \pi t}{n_1^*n_2^*} \bigg)^{ - s}(2s + \alpha  + \beta  - 1)(2s - 2 + \alpha  + \beta ){y^{2s - 3 + \alpha  + \beta }}G(s)\frac{{ds}}{s}. \nonumber  
\end{align}
Using the Mellin transform of $\cos y$, namely \cite{ober}
\[
\mathcal{M}(\cos y,z) := \int_0^\infty y^{z-1} \cos y dy = \Gamma(z) \cos \bigg( \frac{\pi z}{2}\bigg),
\]
valid for $0 < \real(z) < 1$, turns the $y$-integral into
\begin{align}
   &\int_0^\infty  2\cos y\frac{1}{{2\pi i}}\int_{(\frac{5}{4})} \zeta (2s + \alpha  + \beta )\bigg( \frac{2 \pi t}{n_1^*n_2^*} \bigg)^{ - s} \nonumber \\
   &\quad \times (2s + \alpha  + \beta  - 1)(2s - 2 + \alpha  + \beta )y^{2s - 3 + \alpha  + \beta }G(s)\frac{{ds}}{s}dy \nonumber \\
   &= 2\frac{1}{2\pi i}\int_{(\frac{5}{4})}  \zeta (2s + \alpha  + \beta )\bigg( \frac{2 \pi t}{n_1^*n_2^*} \bigg)^{ - s}(2s + \alpha  + \beta  - 1)(2s - 2 + \alpha  + \beta ) \nonumber \\
   &\quad \times \Gamma (2s - 2 + \alpha  + \beta )\cos \bigg( {\frac{{\pi (2s - 2 + \alpha  + \beta )}}{2}} \bigg)G(s)\frac{{ds}}{s} \nonumber \\
   &= 2\frac{1}{{2\pi i}}\int_{(\frac{1}{4})}  \zeta (2s + \alpha  + \beta )\bigg( \frac{2 \pi t}{n_1^*n_2^*} \bigg)^{ - s}\Gamma (2s + \alpha  + \beta )\cos \bigg( {\frac{{\pi (2s - 2 + \alpha  + \beta )}}{2}} \bigg)G(s)\frac{{ds}}{s}. \nonumber  
\end{align}
We have used the duplication formula for the gamma function. We were able to move the path of integration from $\real(s)=\frac{5}{4}$ to $\real(s) = \frac{1}{4}$ without encountering singularities because the simple pole of $\zeta(2s+\alpha+\beta)$ is canceled by the simple zero of the cosine in the integrand at $s=\frac{1-\alpha-\beta}{2}$. Note that if $\alpha,\beta \to 0$, then $\cos(\pi s)$ would still guarantee the lack of a pole \footnote{In \cite[p. 13]{bcr} it is argued, in addition, that $G(\frac{1}{2})=0$ cancels the pole of $\zeta(2s)$ at $s=\frac{1}{2}$.}. Inserting this into the $t$-integral yields
\begin{align}
  \mathcal{A}_{1,0} + O({T^\varepsilon }) &=  - \sum\limits_{{n_1},{n_2}} {} \frac{{{a_{{n_1}}}{{\overline{a}}_{{n_2}}}({n_1},{n_2})^{1-\alpha-\beta}}}{{n_1^{1 - \beta }n_2^{1 - \alpha }}}\int_{ - \infty }^\infty  {} \frac{1}{{2\pi i}}\int_{(\tfrac{1}{4})} {} 2\zeta (2s + \alpha  + \beta )\bigg( \frac{2 \pi t}{n_1^*n_2^*} \bigg)^{ - s} \nonumber \\
   &\quad \times \Gamma (2s + \alpha  + \beta )\cos \bigg( {\frac{{\pi (2s - 2 + \alpha  + \beta )}}{2}} \bigg)G(s)\frac{{ds}}{s}\Phi \bigg( {\frac{t}{T}} \bigg)\frac{{dt}}{{{t^{\alpha  + \beta }}}} \nonumber \\
   &= 2\frac{\pi }{2}\sum\limits_{{n_1},{n_2}} {} \frac{{{a_{{n_1}}}{{\overline{a}}_{{n_2}}}({n_1},{n_2})^{1-\alpha-\beta}}}{{n_1^{1 - \beta }n_2^{1 - \alpha }}}\int_{ - \infty }^\infty  {} \frac{1}{{2\pi i}}\int_{(\tfrac{1}{4})}  {2^{2s + \alpha  + \beta }}{\pi ^{2s + \alpha  + \beta  - 1}}\bigg( \frac{2 \pi t}{n_1^*n_2^*} \bigg)^{ - s} \nonumber \\
   &\quad \times \zeta (1 - (2s + \alpha  + \beta ))G(s)\frac{{ds}}{s}\Phi \bigg( {\frac{t}{T}} \bigg)\frac{{dt}}{{{t^{\alpha  + \beta }}}} \nonumber \\
   &=  - {(2\pi )^{\alpha  + \beta }}\sum\limits_{{n_1},{n_2}}  \frac{{{a_{{n_1}}}{{\overline{a}}_{{n_2}}}({n_1},{n_2})^{1-\alpha-\beta}}}{{n_1^{1 - \beta }n_2^{1 - \alpha }}} \nonumber \\
   &\quad \times \int_{ - \infty }^\infty  \frac{1}{{2\pi i}}\int_{( - \tfrac{1}{4})}  \zeta (1 + 2w - \alpha  - \beta ) \bigg( \frac{t}{2 \pi n_1^*n_2^*} \bigg)^{w} G(w)\frac{{dw}}{w}\Phi \bigg( {\frac{t}{T}} \bigg)\frac{{dt}}{{{t^{\alpha  + \beta }}}}, \nonumber  
\end{align}
after an application of the functional equation of the Riemann zeta-function and the change of variables $s \to -w$ noting that $G(w)=G(-w)$. From the diagonal terms we have
\begin{align} \label{D1A1}
\mathcal{D}_1 + \mathcal{A}_{1,0} &= \sum\limits_{{n_1},{n_2}} {} \frac{{{a_{{n_1}}}{{\overline{a}}_{{n_2}}}{{({n_1},{n_2})}^{1 + \alpha  + \beta }}}}{{n_1^{1 + \alpha }n_2^{1 + \beta }}}\int_{ - \infty }^\infty  {} \frac{1}{{2\pi i}}\int_{(2)} {} {\bigg( {\frac{{2\pi n_1^*n_2^*}}{t}} \bigg)^{ - s}}\zeta (1 + \alpha  + \beta  + 2s)G(s)\frac{{ds}}{s}\Phi \bigg( {\frac{t}{T}} \bigg)dt \nonumber \\
&- {(2\pi )^{\alpha  + \beta }}\sum\limits_{{n_1},{n_2}} {} \frac{{{a_{{n_1}}}{{\overline{a}}_{{n_2}}}({n_1},{n_2})^{1-\alpha-\beta}}}{{n_1^{1 - \beta }n_2^{1 - \alpha }}} \nonumber \\
   &\quad \times \int_{ - \infty }^\infty  {} \frac{1}{{2\pi i}}\int_{( - \tfrac{1}{4})} {} \zeta (1 + 2w - \alpha  - \beta ) \bigg( \frac{t}{2 \pi n_1^*n_2^*} \bigg)^{w}G(w)\frac{{dw}}{w}\Phi \bigg( {\frac{t}{T}} \bigg)\frac{{dt}}{{{t^{\alpha  + \beta }}}} + O(T^\varepsilon).  
\end{align}
Now, we move the path of integration of the $\mathcal{D}_1$ integral to $\real(s)=-\frac{1}{4}$ and we only pick up a simple pole at $s=0$ for which
\[
\mathop {\operatorname{res} }\limits_{s = 0} {x^{ - s}}\zeta (1 + \alpha  + \beta  + 2s)\frac{{G(s)}}{s} = G(0)\zeta (1 + \alpha  + \beta ) = \zeta (1 + \alpha  + \beta )
\]
since $G(0)=1$. This is the only singularity since the pole of $\zeta(1+\alpha+\beta+2s)$ at $s=-\frac{\alpha+\beta}{2}$ is canceled by the simple zero of $G(s)$ at $s=-\frac{\alpha+\beta}{2}$.

We take the chance to explain what happens with $I_2$. Recall that
\[
I_2(\alpha ,\beta ) = \sum_{m_1,m_2,n_1,n_2} \frac{{{a_{{n_1}}}{{\bar a}_{{n_2}}}}}{{m_1^{1/2 - \beta }m_2^{1/2 - \alpha }n_1^{1/2}n_2^{1/2}}}\int_{ - \infty }^\infty  {} {\left( {\frac{{{m_1}{n_2}}}{{{m_2}{n_1}}}} \right)^{it}}{V_{ - \beta , - \alpha }}({m_1}{m_2},t) X_{\alpha,\beta,t} \Phi \bigg(\frac{t}{T}\bigg) dt.
\]
We split into diagonal and off-diagonal cases. In the diagonal case we can immediately use the approximation \eqref{Xdefandasymp}. In the off-diagonal case we truncate the sum by means of the rapid decay of $V_{\alpha,\beta}(x,t)$, then integrate by parts. Here we use the fact that
\begin{align*}
\frac{\partial^j}{\partial t^j} X_{\alpha, \beta,t} &\ll_j t^{-j},
\end{align*}
which follows from Cauchy's integral formula and Stirling's approximation for $\Gamma(z)$. Having done so, we may then use \eqref{Xdefandasymp} and bound the error as we did with the error in the approximation \eqref{gasymp}. A similar analysis then shows that the diagonal terms $\mathcal{D}_2$ of $I_2$ are given by
\begin{align}
\mathcal{D}_2 &= \sum\limits_{{n_1},{n_2}} {} \frac{{{a_{{n_1}}}{{\bar a}_{{n_2}}}{{({n_1},{n_2})}^{1 - \alpha  - \beta }}}}{{n_1^{1 - \beta }n_2^{1 - \alpha }}} \nonumber \\
& \quad \times \int_{ - \infty }^\infty  \frac{1}{{2\pi i}}\int_{(2)} {} {\left( {\frac{{2\pi n_1^*n_2^*}}{t}} \right)^{ - s}}\zeta (1 - \alpha  - \beta  + 2s)G(s)\frac{{ds}}{s} \bigg(\frac{t}{2\pi}\bigg)^{-\alpha-\beta}\Phi \left( {\frac{t}{T}} \right)dt.
\end{align}
The off-diagonal terms coming from $I_2$ are given by
\begin{align}
   \mathcal{A}_{2,0} &= \sum_{d \le N} \sideset{}{'}\sumtwo_{N_1 \le N, N_2 \le N} \sideset{}{'}\sum_{M \le T^{1/2+\varepsilon} \sqrt{\frac{N_2}{N_1}}} \sum_{0 < |\Delta| \le \frac{D}{d}} \sum_{\substack{n_1, n_2 \\ (n_1,n_2)=1}} \frac{a_{dn_1}\bar{a}_{dn_2}F_{N_1}(dn_1)F_{N_2}(dn_2)}{d n_1^{1+\alpha} n_2^{1+\beta}} \nonumber \\
& \quad \times \int_{-\infty}^\infty \int_0^\infty F_M \bigg( \frac{x}{n_1} \bigg) \e \bigg( \frac{\Delta t}{2 \pi x} \bigg) W \left(\frac{2 \pi x^2 }{n_1 n_2 t}\right) \frac{dx}{x^{1-\alpha-\beta}} {\left( {\frac{t}{{2\pi }}} \right)^{ - \alpha  - \beta }}\Phi\bigg( \frac{t}{T}\bigg)dt. \nonumber
\end{align}
Therefore, we see that
\begin{align}
  \mathcal{D}_1 + \mathcal{A}_{2,0} &= \mathcal{M}_1 + O(T^\varepsilon)\nonumber \\
   &\quad + \sum\limits_{{n_1},{n_2}} {} \frac{{{a_{{n_1}}}{{\bar a}_{{n_2}}}{{({n_1},{n_2})}^{1 + \alpha  + \beta }}}}{{n_1^{1 + \alpha }n_2^{1 + \beta }}}\int_{ - \infty }^\infty  {} \frac{1}{{2\pi i}}\int_{( - \tfrac{1}{4})} {} {\left( {\frac{{2\pi n_1^*n_2^*}}{t}} \right)^{ - s}}\zeta (1 + \alpha  + \beta  + 2s)G(s)\frac{{ds}}{s}d\Phi \nonumber \\
   &\quad - \sum\limits_{{n_1},{n_2}} {} \frac{{{a_{{n_1}}}{{\bar a}_{{n_2}}}{{({n_1},{n_2})}^{1 + \alpha  + \beta }}}}{{n_1^{1 + \alpha }n_2^{1 + \beta }}}\int_{ - \infty }^\infty  {} \frac{1}{{2\pi i}}\int_{( - \tfrac{1}{4})} {} {\left( {\frac{{2\pi n_1^*n_2^*}}{t}} \right)^{ - s}}\zeta (1 + \alpha  + \beta  + 2s)G(s)\frac{{ds}}{s}d\Phi  \nonumber \\
   &= \mathcal{M}_1 + O(T^\varepsilon)\nonumber  
\end{align}
where
\[
\mathcal{M}_1 = \sum\limits_{{n_1},{n_2}} {} \frac{{{a_{{n_1}}}{{\bar a}_{{n_2}}}{{({n_1},{n_2})}^{1 + \alpha  + \beta }}}}{{n_1^{1 + \alpha }n_2^{1 + \beta }}}\int_{ - \infty }^\infty  {} \zeta (1 + \alpha  + \beta )\Phi \left( {\frac{t}{T}} \right)dt.
\]
Likewise, $\mathcal{D}_2 + \mathcal{A}_{1,0} = \mathcal{M}_2 + O(T^\varepsilon)$, where
\[
\mathcal{M}_2 = \sum\limits_{{n_1},{n_2}} {} \frac{{{a_{{n_1}}}{{\bar a}_{{n_2}}}{{({n_1},{n_2})}^{1 - \alpha  - \beta }}}}{{n_1^{1 - \beta }n_2^{1 - \alpha }}}\int_{ - \infty }^\infty  {} \zeta (1 - \alpha  - \beta )\Phi \left( {\frac{t}{T}} \right){\left( {\frac{t}{{2\pi }}} \right)^{ - \alpha  - \beta }}dt.
\]
To account for the arithmetical terms in front of the integral we have used the fact that
\[
\frac{{{a_{{n_1}}}{{\overline{a}}_{{n_2}}}{{({n_1},{n_2})}^{1 - \alpha  - \beta }}}}{{n_1^{1 - \beta }n_2^{1 - \alpha }}} = {\bigg( {\frac{{[{n_1},{n_2}]}}{{({n_1},{n_2})}}} \bigg)^{\alpha  + \beta }}\frac{{{a_{{n_1}}}{{\overline{a}}_{{n_2}}}{{({n_1},{n_2})}^{1 + \alpha  + \beta }}}}{{n_1^{1 + \alpha }n_2^{1 + \beta }}},
\]
since $[n_1,n_2](n_1,n_2)=n_1n_2$. Consequently, the total contribution to the main terms coming from the diagonal terms and the two contributing pieces of the off-diagonal terms is
\begin{align}
   \mathcal{D}_1+\mathcal{D}_2 + \mathcal{A}_{1,0} + \mathcal{A}_{2,0} &= \sum\limits_{{n_1},{n_2} \le N}  \frac{{{a_{{n_1}}}{{\overline{a}}_{{n_2}}}{{({n_1},{n_2})}^{1 + \alpha  + \beta }}}}{{n_1^{1 + \alpha }n_2^{1 + \beta }}} \nonumber \\
	 &\quad \times \int_{ - \infty }^\infty  {} \bigg( {\zeta (1 + \alpha  + \beta ) + \zeta (1 - \alpha  - \beta ){{\bigg( {\frac{{2\pi }}{t}\frac{{{n_1}{n_2}}}{{{{({n_1},{n_2})}^2}}}} \bigg)}^{\alpha  + \beta }}} \bigg)\Phi \bigg( {\frac{t}{T}} \bigg)dt \nonumber \\
	 &+O(T^\varepsilon). \nonumber  
\end{align}
This explains the main term of Theorem \ref{theorem1}. In the next sections we estimate the error terms.

\subsubsection{The case $|h| \ge H_d$.}
Recall that $H_d:= \frac{N^2}{dM}T^\varepsilon$. We make the change variables $t=xy$ so that
\begin{align}
  \frac{1}{n_1 n_2}  \tilde{\mathcal{A}}_{M,N_i}(h,{n_1},{n_2},\Delta ) &= \frac{1}{{{n_1}{n_2}}}\int_{ - \infty }^\infty  \operatorname{e} \left( \frac{\Delta y}{2\pi } \right) \nonumber \\
   &\quad \times \int_0^\infty  \operatorname{e} \left( { - \frac{{hx}}{{{n_1}{n_2}}}} \right){F_M}\left( {\frac{x}{{{n_1}}}} \right) W\bigg(\frac{2 \pi x}{n_1 n_2 y}\bigg) \Phi \left( {\frac{{xy}}{T}} \right)\frac{{dx}}{{{x^{\alpha  + \beta }}}}dy .\nonumber 
\end{align} 
Since $F_M$ is supported in $[M/2,3M]$ we have that $x \asymp \frac{N_1M}{d}$. Moreover, $\frac{y}{T} \asymp \frac{1}{x} \asymp \frac{d}{N_1M}$ since $\Phi$ is supported in the interval $[1,2]$. Furthermore, $\frac{1}{n_1 n_2 y} \ll \frac{T^{\varepsilon_1}}{x} \asymp \frac{dT^{\varepsilon_1}}{N_1 M}$ due to the rapid decay of $W$. Integrating by parts $\ell$ times, we obtain
\begin{align}
  &\frac{1}{{{n_1}{n_2}}}\int_0^\infty  {} \operatorname{e} \left( { - \frac{{hx}}{{{n_1}{n_2}}}} \right){F_M}\left( {\frac{x}{{{n_1}}}} \right) W\bigg(\frac{2 \pi x}{n_1 n_2 y}\bigg) \Phi \left( {\frac{{xy}}{T}} \right)\frac{{dx}}{{{x^{\alpha  + \beta }}}} \nonumber \\
  & \ll _{\ell ,\varepsilon }\frac{{{d^2}}}{{{N_1}{N_2}}}{\left( {\frac{{{n_1}{n_2}}}{h}\frac{{d{T^{{\varepsilon _1}}}}}{{M{N_1}}}} \right)^{\ell  + 1}}\frac{{M{N_1}}}{d} \ll {\left( {\frac{{{T^{{\varepsilon _1}}}}}{h}} \right)^{\ell  + 1}}\left( {\frac{{{N_2}}}{{dM}}} \right)^\ell . \nonumber 
\end{align}
Therefore, the contribution to $\mathcal{S}_1$ from $|h|>H_d$ is
\begin{align}
& \ll \sum_{d \le N} \sideset{}{'}\sum\limits_{\substack{N_1, N_2 \le N \\ M \le T^{1/2+\varepsilon} \sqrt{\frac{N_2}{N_1}}}} \sum_{0 \le |\Delta| \le \frac{D}{d}} \sum_{\substack{n_1, n_2 \\ (n_1, n_2)=1}} \frac{a_{dn_1}\bar{a}_{dn_2}F_{N_1}(dn_1)F_{N_2}(dn_2)}{d n_1^{1-\beta} n_2^{1-\alpha}} \sum_{|h| \ge H_d} \frac{dT}{hN_1M} \bigg( \frac{N_2T^{\varepsilon_1}}{dMh} \bigg)^\ell  \ll T^{-A}, \nonumber 
\end{align}
when $\ell$ is sufficiently large. Thus, the terms for which $|h|>H_d$ yield a negligible contribution.

\subsubsection{The case $0< |h| < H_d$.}
It is sufficient to consider the terms $0 < h < H_d$. We change of variables $t=xy$, followed by $x \to x n_1 n_2$, and consider 
\begin{align}
  \mathcal{A}_{M,N_1,N_2}^* &:= \sum_{\substack {n_1,n_2 \\ (n_1,n_2) = 1}} \sum_{0 < |\Delta | < \frac{D}{d}}  \sum_{0 < h < H_d}  \frac{a_{dn_1}\overline{a}_{dn_2}F_{N_1}(dn_1)F_{N_2}(dn_2)}{dn_1^\alpha n_2^\beta }\operatorname{e} \left(  - \frac{h\Delta \bar n_1}{n_2} \right) \nonumber \\
   &\quad \times \int_0^\infty   \operatorname{e} ( - hx){F_M}(xn_2)\int_{ - \infty }^\infty   \operatorname{e} \left( \frac{\Delta y}{2\pi} \right) W\bigg( \frac{2 \pi x}{y} \bigg)\Phi \left( \frac{xyn_1n_2}{T} \right)dy\frac{dx}{x^{\alpha  + \beta }}. \nonumber  
\end{align}
To decouple the variables $n_1$ and $n_2$ we write $\Phi$ in terms of its Mellin transform $\mathcal{M}(\Phi,w)$, i.e. 
\[
\Phi \left( \frac{xyn_1n_2}{T} \right) = \frac{1}{2\pi i}\int_{(\varepsilon )}  \mathcal{M}(\Phi ,w)\left( \frac{xyn_1n_2}{T} \right)^{ - w}dw.
\]
Let $h \Delta = a$, $A=\frac{DH_d}{d}=\frac{N_1N_2}{d^2T^{1-\varepsilon}}$, and $\nu_{x,y}(a) = \sum_{h\Delta=a} \operatorname{e} (-hx + \frac{\Delta y}{2 \pi})$. With this notation, we arrive at the following
\begin{align}
  \mathcal{A}_{M,N_1,N_2}^* &= \frac{1}{2\pi i d}\int_0^\infty   \int_{ - \infty }^\infty   \int_{(\varepsilon )}  W\bigg(  \frac{2 \pi x}{y}\bigg)\sum\limits_{0 < |a| < A}  \nu _{x,y}(a) \nonumber \\
   &\quad \times \sum_{\substack{n_1,n_2 \\ (n_1,n_2) = 1}}  \frac{a_{dn_1}\overline{a}_{dn_2}F_{N_1}(dn_1)F_{N_2}(dn_2)F_M(x{n_2})}{n_1^{\alpha  + w}n_2^{\beta  + w}}\operatorname{e} \left(  - \frac{a\bar n_1}{n_2} \right)\mathcal{M}(\Phi ,w)\frac{T^w}{x^wy^w}dwdy\frac{dx}{x^{\alpha  + \beta }}. \nonumber  
\end{align} 
Observe that since $F_M$ is supported in $[M/2,3M]$ we have $x \asymp \frac{dM}{N_2}$. Moreover, $y \asymp \frac{T}{xn_1n_2} \asymp \frac{Td}{MN_1}$ because $\Phi$ is supported in $[1,2]$. We now distinguish three cases.

The first, and easiest, case is when we have no information about the coefficients $a_n$ other than $a_n \ll_\varepsilon n^\varepsilon$. Here we use Lemma \ref{BClemma}, as in \cite{bcr}. In our slightly modified setting we have only to note that $T^{-O(1)} < x,y < T^{O(1)}$ and $|\alpha|,|\beta| \ll \frac{1}{\log T}$. The second and third cases, in which we specialize the coefficients, are more difficult and we give the proofs after we bound the error $\mathcal{E}_W$.

\subsubsection{The bound for $\mathcal{E}_W$}
Here we bound the error $\mathcal{E}_W$, which appeared in \eqref{definitionofE}. As with $\mathcal{A}_1$ we extract the common divisor $d$ from $n_1$ and $n_2$, apply the Poisson summation formula, and change variables to obtain
\[
\mathcal{E}_W = \sum_{d \le N} \frac{1}{d} \sideset{}{'}\sum\limits_{N_1, N_2 \le N}  \sideset{}{'}\sum\limits_{M \le T^{1/2+\varepsilon} \sqrt{\frac{N_2}{N_1}}} \sum_{0 < |\Delta| \le \frac{D}{d}} \sum_{\substack{n_1, n_2 \\ (n_1,n_2)=1}} a_{dn_1}\bar{a}_{dn_2} F_{N_1}(dn_1)F_{N_2}(dn_2) \mathcal{E}_{M,N_i}(n_1,n_2,\Delta), 
\]
where
\begin{align}
\mathcal{E}_{W,M,N_i}(n_1,n_2,\Delta) &= \frac{1}{n_1^{1-\beta}n_2^{1-\alpha}} \sum_{h \in \mathbb{Z}} \operatorname{e}\bigg(-\frac{h\bar{n}_1 \Delta}{n_2}\bigg) \int_0^\infty \operatorname{e}\bigg(-\frac{hx}{n_1n_2}\bigg) F_M \bigg(\frac{x}{n_1}\bigg) \int_{-\infty}^{\infty} \operatorname{e}\bigg(\frac{\Delta t}{2 \pi x}\bigg) \Phi \bigg(\frac{t}{T}\bigg) \nonumber \\
&\quad \times \bigg[ W\bigg(\frac{2 \pi x^2}{n_1n_2t}\bigg)\bigg(-\frac{(\tfrac{1}{2}+\alpha)\Delta}{x^2}-\frac{it\Delta^2}{2x^3}\bigg) + \frac{2\pi \Delta}{n_1n_2t}W'\bigg(\frac{2 \pi x^2}{n_1n_2t}\bigg)\bigg]dt\frac{dx}{x^{\alpha+\beta}}. \nonumber
\end{align}
If we integrate by parts, as in Case 2 of $\mathcal{A}_1$, we see that the contribution coming from the terms with $|h| >H_d$ is $O(1)$. The rest of the proof is finished by estimating trivially the remaining terms, namely
\[
\mathcal{E}_{W,M,N_i}(n_1,n_2,\Delta) \ll \frac{T^\varepsilon}{n_1 n_2} \bigg(1+\frac{N_2}{dM}\bigg),
\]
from which we obtain
\[
\mathcal{E}_W \ll T^{-1/2+\varepsilon}N + T^{-1+\varepsilon}N^2 \ll NT^{\varepsilon}.
\]

\subsubsection{Bounding the Error Terms: Feng and Conrey}\label{confeng error terms}
For specialized coefficients $a_n$, we study here the sum
\begin{align*}
\mathcal{B} = \mathcal{B}(A,N_1,N_2,d) &= \sum_{0 < |a| < A} \nu(a) \sum_{\substack{n_1 \leq N \\ n_2 \asymp N_2/d \\ (n_1,n_2) = 1}} \frac{a_{dn_1} F_{N_1}(dn_1)}{n_1^{\alpha + w}}r(n_2) e \left(-a \frac{\overline{n_1}}{n_2} \right),
\end{align*}
where $\nu,r$ are functions satisfying $\nu(n),r(n) \ll_\varepsilon n^\varepsilon$. We give full details only when the coefficients $a_n$ are coefficients of the Feng mollifier. The argument is virtually identical in the case of the Conrey mollifier, and we indicate some of these differences as we go along. Our argument is based on that of Conrey \cite{conrey89}.

Let us now suppose that the $a_n$ are given by the coefficients of the Feng mollifier (see, e.g. \cite{feng}), that is,
\begin{align*}
a_n &= \sum_{2 \leq k \leq K} \frac{1}{(\log N)^k}\mu^2(n)(\mu * \Lambda^{*k})(n) P_k\left(\frac{\log(N/n)}{\log N} \right),
\end{align*}
where the $P_k$ are polynomials satisfying certain properties and $K$ is a fixed integer. By linearity we see it suffices to study
\begin{align*}
a_n = \mu^2(n)(\mu * \Lambda^{*k})(n) P_k\left(\frac{\log(N/n)}{\log N} \right).
\end{align*}

What shows up in $\mathcal{B}$ is not $a_{n_1}$, but $a_{dn_1}$, and we need to separate $d$ and $n_1$ from one another as much as possible. It is easy to separate $d$ and $n_1$ inside of $P_k$: by linearity and the binomial theorem we reduce to studying
\begin{align*}
a_{dn_1} &= \mu^2(dn_1)(\mu * \Lambda^{*k})(dn_1) (\log n_1)^j,
\end{align*}
for some integers $j,k \geq 0$. The presence of the $\mu^2$ factor means we may assume $(d,n_1) = 1$, and thus $\mu^2(dn_1) = \mu^2(d) \mu^2(n_1)$. It therefore remains to separate $d$ and $n_1$ in $(\mu * \Lambda^{*k})(dn_1)$.

For coprime integers $u,g$, we have
\begin{align*}
(\mu* \Lambda^{*j})(ug) &= \mathop{\sum \cdots \sum}_{n \ell_1 \cdots \ell_j = ug} \mu(n) \Lambda(\ell_1) \cdots \Lambda(\ell_j).
\end{align*}
Since $(u,g) = 1$ we have $\ell_i | u$ or $\ell_j | g$, but we cannot have $\ell_j |u$ and $\ell_j | g$. It follows that $(\mu* \Lambda^{*j})(ug)$ is the sum of $2^j$ sums of the form
\begin{align*}
\mathop{\sum \cdots \sum}_{\substack{n \ell_1 \cdots \ell_j = ug \\ \ell_{i_1}, \cdots, \ell_{i_s} | u \\ \ell_{i_{s+1}}, \ldots, \ell_{i_j} | g}} \mu(n) \Lambda(\ell_1) \cdots \Lambda(\ell_j).
\end{align*}
Since $n | ug$ with $(u,g) = 1$ we may uniquely write $n = n_u n_g$, where $n_u | u$ and $n_g | g$. Obviously $(n_u,n_g) = 1$. Moreover, it is easy to see that $n_u \ell_{i_1} \cdots \ell_{i_r} = u$. We therefore have
\begin{align*}
\mathop{\sum \cdots \sum}_{\substack{n \ell_1 \cdots \ell_j = ug \\ \ell_{i_1}, \cdots, \ell_{i_s} | u \\ \ell_{i_{s+1}}, \ldots, \ell_{i_j} | g}} \mu(n) \Lambda(\ell_1) \cdots \Lambda(\ell_j) &= (\mu * \Lambda^{*s})(u) (\mu * \Lambda^{*(j-s)})(g),
\end{align*}
and this gives the desired separation of $u$ and $g$.

It follows that $\mathcal{B}$ is a linear combination of $O(1)$ sums of the form
\begin{align*}
\tilde{\mathcal{B}} &= \lambda(d)\sum_{0 < |a| < A} \nu(a) \sum_{\substack{n_1 \leq N \\ n_2 \asymp N_2/d \\ (n_1,dn_2) = 1}} \frac{\mu^2(n_1) (\mu*\Lambda^{*s})(n_1) F_{N_1}(dn_1)(\log n_1)^j}{n_1^{\alpha + w}}r(n_2) e \left(-a \frac{\overline{n_1}}{n_2} \right),
\end{align*}
where $\lambda(d) \ll (dT)^\varepsilon$ and $j,s \geq 0$ are fixed integers. Observe that $\mu * \Lambda^{*0} = \mu$.

Before proceeding, it is helpful to slightly clean up the notation. We set $U = N_1/d$ and $V = N_2/d$, so that we need to estimate
\begin{align*}
\mathcal{B}_1 &= \sum_{0 < |a| < A} \nu(a) \sum_{\substack{u \leq N \\ u \asymp U \\ v \asymp V \\ (u,dv) = 1}} \frac{\mu^2(u) (\mu*\Lambda^{*s})(u) F_{dU}(du)(\log u)^j}{u^{\alpha + w}}r(v) e \left(-a \frac{\overline{u}}{v} \right).
\end{align*}
The next step is to decompose $\mu*\Lambda^{*s}$ into different pieces. This will give rise to Type I and Type II sums, as they are often called in the literature. We recall the following identities, due essentially to Heath-Brown \cite{heathbrown}, for $\mu$ and $\Lambda$, valid for $n \leq 2U$:
\begin{align*}
\Lambda(n) &= \sum_{1 \leq k \leq K} (-1)^{k-1} {K \choose k} \mathop{\sum \cdots \sum}_{\substack{m_1 \cdots m_k n_1 \cdots n_k=n \\ m_1,\ldots,m_k \leq (2U)^{1/K}}} \mu(m_1) \cdots \mu(m_k) \log(n_k), \\
\mu(n) &= \sum_{1 \leq k \leq K} (-1)^{k-1} {K \choose k} \mathop{\sum \cdots \sum}_{\substack{m_1 \cdots m_k n_1 \cdots n_{k-1}=n \\ m_1,\ldots,m_k \leq (2U)^{1/K}}} \mu(m_1) \cdots \mu(m_k).
\end{align*}
We apply these identities with $K=2$. We split the range of summation of each variable $m_i,n_j$ into dyadic intervals of the form $X < x \leq 2X$, which implies that for $U <u \leq 2U$ the function $(\mu* \Lambda^{*s})(u)$ is a linear combination of $O((\log U)^{4s+3})$ functions of the form
\begin{align*}
\mathop{\sum \cdots \sum}_{\substack{n_1 \cdots n_{4s+3} =u \\ n_i \in I_i}} \mu(n_1) \cdots \mu(n_{2s+2}) \log(n_{2s+3}) \cdots \log(n_{3s+2}).
\end{align*}
Here $I_i = (X_i,2X_i]$, $2^{-(4s+3)}U \leq \prod_i X_i < 2U$, and $2X_i \leq (2U)^{1/2}$ for $1 \leq i \leq 2s+2$. It is possible that some $I_i$ contain only the integer 1. 

Let $1 \leq W \ll U^{1/3}$ be a parameter to be chosen. We claim that either there is some $i \in \{1,\ldots,4s+3\}$ with $X_i \gg U/W$, or there is a subset $S \subset \{1,\ldots,4s+3\}$ such that $ W \ll\prod_{i \in S} X_i \ll U/W$. If there exists an $i$ with $X_i \gg U/W$ we are done, and if there is some $i$ such that $W \ll X_i \ll U/W$ we are also done (take $S = \{i\}$). Thus we may suppose $X_i \ll W$ for all $i$. Since $X_i \ll W$ and $\prod_i X_i \gg U \gg W$, there is some minimal $i_0 \geq 2$ such that
\begin{align*}
\prod_{i = 1}^{i_0} X_i \gg W.
\end{align*}
By minimality we have
\begin{align*}
\prod_{i = 1}^{i_0} X_i &= X_{i_0}\prod_{i = 1}^{i_0-1} X_i \ll X_{i_0} W \ll W^2 \ll \frac{U}{W},
\end{align*}
the last inequality following since $W \ll U^{1/3}$. We finish by taking $S = \{1,\ldots,i_0\}$. To balance the various error terms arising we eventually take $W = U^{1/6}$.

It follows that $(\mu*\Lambda^{*s})(u)$ is a linear combination of $O(U^\varepsilon)$ functions of the form $(\beta * g) (u)$, where $\beta$ is supported on integers $\ll W$ and $g$ is equal to the constant one function or $\log$ (the Type I case), or functions of the form $\gamma * \delta$, where $\gamma,\delta$ are supported on integers $W \ll n \ll U/W$ (the Type II case). The functions $\beta,g,\gamma,\delta$ are supported on dyadic intervals, and satisfy the bounds $\beta(n),g(n),\gamma(n),\delta(n) \ll_\varepsilon n^\varepsilon$.

In dealing with the Conrey mollifier we perform a similar combinatorial decomposition on the M\"obius function, and introduce a similar parameter $W'$, which is eventually taken to be $W' = U^{1/4}$.

Let us first consider a Type I sum. Using the binomial theorem to separate variables inside the logarithm, we must therefore estimate
\begin{align*}
\sum_{a \leq A} \nu(a) \sum_{v \asymp V} r(v)\sum_{\substack{e \asymp E \\ (e,dv) = 1}} \beta(e) \sum_{\substack{f \asymp F \\ ef \asymp U \\ ef \leq N \\ (f,dev) = 1}} \mu^2(f) \frac{(\log f)^\ell F_{dU}(def)}{f^{\alpha + w}} e \left(-a \frac{\overline{ef}}{v} \right),
\end{align*}
where $EF \asymp U$, $E \ll W$, and $\ell \geq 0$ is an integer. By summation by parts, we have
\begin{align*}
\sum_{\substack{f \asymp F \\ ef \asymp U \\ ef \leq N \\ (f,dev) = 1}} \mu^2(f) \frac{(\log f)^\ell F_{dU}(def)}{f^{\alpha + w}} e \left(-n \frac{\overline{ef}}{v} \right) &\ll (1+|w|) T^\varepsilon \sum_{\substack{f \in I \\ f \asymp F \\ (f,dev)=1}} \mu^2(f) e \left(-n \frac{\overline{ef}}{v} \right)
\end{align*}
for some interval $I$. By inclusion-exclusion this latter sum is equal to
\begin{align*}
\sum_{\substack{h \ll F^{1/2} \\ (h,dev) = 1}}\mu(h) \sum_{\substack{f_1 \in I/h^2 \\ f_1 \asymp F/h^2 \\ (f_1,dev)=1}} e \bigg(-n \frac{\overline{ef_1 h^2}}{v} \bigg).
\end{align*}
The inner sum is trivially $\ll F/h^2$. By Weil's bound for Kloosterman sums \cite{weilbound}, the inner sum is also $\ll T^\varepsilon v^{1/2} (n,v)(1 + F v^{-1})$. Taking the minimum of these two bounds and using the inequality $\min(x,y) \leq (xy)^{1/2}$, we have
\begin{align*}
\sum_{\substack{f \in I \\ f \asymp F \\ (f,dev)=1}} \mu^2(f) e \left(-n \frac{\overline{ef}}{v} \right) &\ll F^{1/2} v^{1/4} (1 + F^{1/2} v^{-1/2} ) (n,v)^{1/2} T^\varepsilon.
\end{align*}
On summing over $a,v$, and $e$ we obtain that the contribution to $\mathcal{B}_1$ from a Type I sum is
\begin{align*}
\ll (1+|w|)T^\varepsilon A (W^{1/2} U^{1/2} V^{5/4} + UV^{3/4} ).
\end{align*}

In the case of Conrey's mollifier we also arrive at incomplete Kloosterman sums, but now the summation variable is not weighted by a factor $\mu^2(f)$. We are therefore able to apply Weil's bound for Kloosterman sums directly.

We turn now to studying Type II sums. Separating variables via the Mellin transform of $F$ and the binomial theorem, it suffices to bound the sum
\begin{align*}
\sum_{v \asymp V} r(v) \sum_{a \leq A} \nu(a) \sum_{\substack{\substack{b \asymp B \\ (b,v)=1}}} \gamma(b) \sum_{\substack{c \asymp C \\ (c,bv)=1}} \delta(c) e\bigg(a \frac{\overline{bc}}{v} \bigg),
\end{align*}
say, where $BC \asymp U$, and $W \ll B,C \ll U/W$. We may assume without loss of generality that $B \ll U^{1/2}$, so that in fact $W \ll B \ll U^{1/2}$. This is almost in a form where we may apply Lemma \ref{DIlemma}, but we have the condition $(b,c) = 1$. However, this condition may be removed with M\"obius inversion at no cost. We deduce that the contribution to $\mathcal{B}_1$ from a Type II sum is
\begin{align*}
\ll T^\varepsilon \bigg(\sum_{(b,a,u,v) \in J} B^b A^a U^u V^v \bigg)^{1/4},
\end{align*}
where
\begin{align*}
J &= \{(-2,2,4,4),(-1,4,3,4),(1,4,3,3),(1,2,3,4),(3,2,3,3),(0,4,4,2),(1,3,4,2)\}.
\end{align*}
(The non-alphabetic ordering of the components of the tuples of $J$ is to facilitate comparison with \cite[p. 23]{conrey89}). With Conrey's mollifier we bound the Type II sums in the same fashion, but here there is no need for M\"obius inversion to remove a coprimality condition.

Combining our bounds and integrating over $x,y,w$, we find that
\begin{align*}
\mathcal{A}_{M,N_1,N_2}^* &\ll \frac{T^\varepsilon}{dA} \bigg(A W^{1/2} U^{1/2} V^{5/4} + AUV^{3/4} +  \bigg(\sum_{(b,a,u,v) \in J} B^b A^a U^u V^v \bigg)^{1/4} \bigg).
\end{align*}
We set $W = U^{1/6}$ and recall that $W \ll B \ll U^{1/2}, U = N_1/d, V = N_2/d$. Summing over dyadic intervals $M \ll T^{O(1)}$, $N_i \leq N$, and $d \leq N$, we find that the contribution to $\mathcal{A}$ from $0 < |h| \leq H_d$ is bounded by
\begin{align*}
&\ll T^\varepsilon (N^{11/6} + N^{11/12}T^{1/2} ).
\end{align*}
If $N = T^\theta$ this error is $\ll T^{1-\varepsilon}$ for $\theta < \frac{6}{11}$. For Conrey's mollifier one finds that the error term is
\begin{align*}
\ll T^\varepsilon ( N^{7/4} + T^{1/2}N^{7/8} ),
\end{align*}
which is sufficiently small provided $N = T^\theta$ with $\theta < \frac{4}{7}$.


\section{Proof Theorem \ref{theorem2} and Theorem \ref{theorem3}}
\subsection{Proof of Theorem \ref{theorem2}}
The strategy of the main terms is identical to the proof of the previous case with $\overline{a}_{n_2}$ replaced by $\overline{b}_{n_2}$. The difference is in the error term involving the case $0 < |h| < H_d$ when the structure of $A$ and $B$ are different. Recalling that
\begin{align}
  \mathcal{A}_{M,N_1,N_2}^* &= \sum_{\substack {n_1,n_2 \\ (n_1,n_2) = 1}} \sum_{0 < |\Delta | < \frac{D}{d}}  \sum_{0 < h < H_d}  \frac{a_{dn_1}\overline{b}_{dn_2}F_{N_1}(dn_1)F_{N_2}(dn_2)}{dn_1^\alpha n_2^\beta } \nonumber \\
   &\quad \times \operatorname{e} \left(  - \frac{h\Delta \bar n_1}{n_2} \right)\int_0^\infty   \operatorname{e} ( - hx){F_M}(xn_2)\int_{ - \infty }^\infty   \operatorname{e} \left( \frac{\Delta y}{2\pi} \right)W\left( \frac{2\pi x}{y} \right)\Phi \left( \frac{xyn_1n_2}{T} \right)dy\frac{dx}{x^{\alpha  + \beta }}. \nonumber  
\end{align}
We finish by a very similar analysis to that of $\mathsection$3.2.5.
\subsection{Proof of Theorem \ref{theorem3}}
Once again the technique is the same, we first apply the approximate functional equation and separate $J$ into $J_1$ and $J_2$. This time we bear in mind the convolution
\[
\mathfrak{a}_d := \sum_{n_1 k_1 = d}a_{n_1} b_{k_1}
\]
when computing the $J_1$ integral
\[
J_1(\alpha ,\beta ) := \sum\limits_{m_1,m_2,n_1,n_2,k_1,k_1}  \frac{a_{n_1}\overline{a}_{n_2}b_{k_1}\overline{b}_{k_2}}{m_1^{1/2 + \alpha }m_2^{1/2 + \beta }(n_1n_2k_1k_2)^{1/2}}\int_{ - \infty }^\infty   \left( \frac{m_1n_2k_2}{m_2n_1k_1} \right)^{it}V_{\alpha ,\beta }(m_1m_2,t)\Phi \left( \frac{t}{T} \right)dt.
\]
A similar analysis to that of \cite[$\mathsection$4]{bcr} ends the proof. We summarize the steps. By following a path like that of the proof of Theorem \ref{theorem1}, we arrive at
\begin{align}
\mathcal{A}^*_{M,N_1,N_1} &= \sum_{\substack{n_1,n_2 \\ (n_1,n_2)=1}} \sum_{0 < |\Delta| < \frac{D}{d}} \sum_{0 < h < H_d} \frac{\mathfrak{a}_{dn_1}\overline{\mathfrak{a}}_{dn_2}F_{N_1}(dn_1)F_{N_2}(dn_2)}{d n^\alpha n^\beta} \nonumber \\
& \quad \times \e \bigg(-\frac{h \Delta \overline{n}_1}{n_2}\bigg) \int_0^\infty \e (-hx) F_M(xn_2) \int_{-\infty}^\infty \e \bigg(\frac{\Delta y}{2 \pi}\bigg) W\bigg(\frac{2 \pi x}{y}\bigg) \Phi \bigg(\frac{x y n_1 n_2}{T}\bigg) dy \frac{dx}{x^{\alpha+\beta}}. \nonumber
\end{align}
Let us now write $\mathfrak{a}_{dn_1}$ as $a_{\mu h}b_{\nu r}$, where $\mu |d^\infty, (d,j)=1, n_1 = \rho r j, \nu = \frac{d}{(\mu,d)}$ and $\rho = \frac{\mu}{(\mu,d)}$. We leave $\mathfrak{a}_{dn_2}$ unchanged. This implies that the quantity we need to bound is
\begin{align}
&\sum_{d \le N}\frac{1}{d} \sideset{}{'}\sum_{N_1,N_2,M} \sum_{\substack{\mu | d^\infty \\ (\nu = d/(\mu,d))}} \sum_{0 < |\Delta| < \frac{D}{d}} \sum_{0 < |h| < H_d} \sum_{(n_2, \rho)=1} \mathfrak{a}_{dn_2} \sum_{(j,dn_2)=1} a_{\mu j} \sum_{(r,n_2)=1} b_{\nu r} \e \bigg(-\frac{h \Delta \overline{\rho r j}}{n_2}\bigg) \nonumber \\
& \quad \times \frac{F_{N_1}(d \rho r j)F_{N_2}(dn_2)}{n_1^\alpha n_2^\beta} \int_0^\infty \e (-hx) F_M(xn_2) \int_{-\infty}^\infty \e \bigg(\frac{\Delta y}{2 \pi}\bigg) W\bigg(\frac{2 \pi x}{y}\bigg) \Phi \bigg(\frac{x y \rho r j n_2}{T}\bigg) dy \frac{dx}{x^{\alpha+\beta}}.\nonumber	
\end{align}
Here the sums over $N_1,N_2$ and $M$ are dyadic sums up to $NK, NK$ and $T^{1/2+\varepsilon}\sqrt{N_2/N_1}$, respectively.

The key instruments are now a separation of the variables $n_2, r, j$ like the one performed in $\mathsection$3.2.5 via the Mellin transforms of the functions $F_{N_1}$ and $\Phi$, followed by an application of Lemma \ref{DIlemma}. For convenience to the reader we remark the following identification of indices: the sum over $\ell$ in the lemma is the sum over $n_2$ with $L = \frac{N_2}{d}$ but the sum over $j$ remains the same with $J \le \frac{N}{\mu}$; moreover, the sum over $u$ is the sum over $h\Delta$ with $U=\frac{N_1N_2}{d^2 T^{1-\varepsilon}}$, and the sum over $v$ becomes the sum over $r$ with $V \le \frac{K}{\nu}$. Lastly, $JV \le \frac{N_1}{d \rho}$. Once the dyadic sum over $M$ is performed, the result of Lemma \ref{DIlemma} implies that the above expression is bounded by 
\[
T^\varepsilon (T^{1/2}N^{3/4}K + T^{1/2}NK^{1/2} + N^{7/4}K^{3/2}),
\]
see \cite[p. 17]{bcr} for further details.


\section{Application to critical zeros}

We mentioned in $\mathsection$1 that one needed $I(\alpha,\beta)$ rather than $I$ in order to compute the percentage of zeros on the critical line. More precisely, let $N(T)$ and $N_0(T)$ be the number of zeros inside the rectangle $0 < \real(s) < 1$ and on the critical line, respectively, both up to height $0 < t < T$, (see e.g. \cite[$\mathsection$9 and $\mathsection$10]{titchmarsh}). The proportion of zeros\footnote{A history of the values of $\kappa$ is documented in \cite{krz01}. Bui, Conrey and Young \cite{bcy} were able to get $41.05\%$. Feng \cite{feng} claimed a value of $41.27\%$, though this was contested in \cite{bui, krz01,rrz01}, and reduced to $41.07\%$ due to an incomplete claim on the error terms. The calculation in this note shows that the length $\theta$ of Feng's mollifier may indeed be taken to be larger than $\frac{1}{2}$,  but pushing $\theta$ past $\frac{6}{11}$ to (perhaps) $\frac{4}{7}$ will require more effort.} on the line is defined as
\[
\kappa := \mathop {\lim \inf }_{T \to \infty} \frac{N_0(T)}{N(T)}.
\]
Littlewood's lemma yields the useful inequality \cite[p. 290]{titchmarsh} and \cite[p. 7]{conrey89}
\begin{align} \label{kappaineq}
\kappa \ge 1 - \frac{1}{R} \log \bigg(\frac{1}{T}\int_1^T |V A (\sigma_0 + it)|^2 dt \bigg) + o(1),
\end{align}
thereby linking the percentage to twisted second moments. Here $V$ is defined by
\[
V(s) := Q \bigg( -\frac{1}{L} \frac{d}{ds} \bigg)\zeta(s),
\]
where $Q(x)$ is a real polynomial satisfying $Q(0)=1$ and $Q(x)+Q(1-x) = \operatorname{constant}$, and $\sigma_0 = \frac{1}{2}-\frac{R}{L}$ (recall that $R$ is a bounded constant of our choice). In this case the Dirichlet polynomial $A(s)$ is chosen to mimic $\zeta(s)^{-1}$ or $(\zeta(s)+\frac{\zeta'(s)}{L})^{-1}$. Rather than computing the integral in \eqref{kappaineq}, it is more useful to compute $\mathfrak{I}$ defined by 
\[
\mathfrak{I} := Q \bigg(\frac{-1}{\log T} \frac{d}{d\alpha}\bigg)Q \bigg(\frac{-1}{\log T} \frac{d}{d\beta}\bigg)I(\alpha,\beta) \bigg|_{\alpha=\beta=-R/L}.
\]

We use a two-piece mollifier $\psi(s) = \psi_1(s) + \psi_2(s)$. We take $\psi_1(s)$ to be the Conrey mollifier, with coefficients given by
\begin{align*}
a_n &= \frac{\mu(n)}{n^{1/2-\sigma_0}} P_1 \left( \frac{\log (N_1/n)}{\log N_1}\right),
\end{align*}
where $P_1$ is a polynomial satisfying some minor conditions, and $N_1 = T^{\theta_1}$ with $\theta_1 < \frac{4}{7}$. We take $\psi_2(s)$ to be the Feng mollifier, with coefficients
\begin{align*}
a_n &= \frac{1}{n^{1/2-\sigma_0}}\sum_{2 \leq k \leq K} \frac{1}{(\log T)^k}\mu^2(n)(\mu * \Lambda^{*k})(n) P_k\left(\frac{\log(N_2/n)}{\log N} \right).
\end{align*}
The polynomials $P_k$ also satisfy some minor conditions, and we are free to choose the integer parameter $K$. We have $N_2 = T^{\theta_2}$ with $\theta_2 < \frac{6}{11}$.

We next open the square in \eqref{kappaineq} and employ Theorem \ref{theorem1} for integrals involving $\psi_i \overline{\psi_i}$ for $i \in \{1,2\}$, and Theorem \ref{theorem2} for the integrals involving $\psi_1 \overline{\psi_2}$ and $\psi_2 \overline{\psi_1}$. The error terms associated with this process will hold uniformly by Cauchy's integral formula \cite[p. 41]{bcy}. For the sums over $d$ and $e$ one first uses the fact that
\[
(d,e)^{1+\alpha+\beta} = \sum_{\substack{h|d \\ h|e}} \sum_{k|h} \mu(h) \bigg(\frac{h}{k}\bigg)^{1+\alpha+\beta} = \sum_{\substack{h|d \\ h|e}} h^{1+\alpha+\beta} F(h,1+\alpha+\beta) \quad \textnormal{with} \quad F(h,s) := \prod_{p|h}(1-p^{-s}),
\]
and then follows the technique of the main term computations given in \cite[$\mathsection$6]{conrey83a} and \cite[p. 13]{conrey89} for Conrey's mollifier and in \cite[$\mathsection$3]{feng} for Feng's mollifier.

We utilize these main term computations in conjunction with the following choice of parameters: take $\theta_1 = \frac{4}{7} - \varepsilon$, $\theta_2=\frac{6}{11} - \varepsilon$, $R=1.3025$ and $K=3$ in the main terms of Feng mollifier (see e.g. \cite[Theorem 2]{feng} \cite[Theorems 1.1, 1.2 and 1.3]{krz01}) as well as
\begin{align}
P_1(x) &= x + 0.327608 x(1-x) - 1.62086 x(1-x)^2 - 0.160377 x(1-x)^3 + 1.29018 x(1-x)^4, \nonumber \\
P_2(x) &= 0.197567x + 2.40831 x^2, \nonumber \\
P_3(x) &= 0.649142x + 1.042x^2, \nonumber \\
Q(x) &= 0.491203 + 0.630413(1-2x) -0.149615 (1-2x)^3 +0.0279997 (1-2x)^5. \nonumber
\end{align}
This yields $\kappa \ge .41491637$. 

\section{Acknowledgments}
 The authors would like to acknowledge Roger Baker, Maksym Radziwi\l{}\l{}, Arindam Roy, and Alexandru Zaharescu for fruitful and helpful comments. The second author would like to thank Maksym Radziwi\l{}\l{} for his hospitality at McGill University.
 
 The authors thank Roger Heath-Brown for discovering a misprint in the statement of Theorem 1.2 in earlier versions of this manuscript.

\end{document}